\documentclass[11pt,a4paper]{article}

\usepackage{mathtools} 
\usepackage{amssymb,amsthm}
\usepackage{bm}
\usepackage{fancyhdr}
\usepackage{xcolor}
\usepackage{keyval}
\usepackage{hyperref}
\usepackage{tensor}
\usepackage{tikz-cd}
\usepackage{txfonts}
\usepackage[textsize=tiny]{todonotes}
\usepackage{enumerate}
\usepackage{almacros}
\usepackage{multicol}
\usepackage{keyval}

\newcommand{\mrank}{\mathrm{mrank}}

\usepackage[left=2cm,right=2cm,top=2cm,bottom=2cm]{geometry}

\makeatletter \define@key{todonotes}{David}[]{\setkeys{todonotes}{color=blue!20, noline, prepend, caption=David}} \makeatother

\makeatletter \define@key{todonotes}{al}[]{\setkeys{todonotes}{color=green!20, noline, prepend, caption=Fons}} \makeatother

\makeatletter \define@key{todonotes}{leo}[]{\setkeys{todonotes}{color=red!20, noline, prepend, caption=Leo}} \makeatother

\ExplSyntaxOn

\newcommand{\cat}[1]{
  \mathcal{\tl_range:nnn {#1} {1} {1}}~
  \mathrm{\tl_range:nnn {#1} {2} {-1}}
}

\newcommand{\cattwo}[1]{
  \mathcal{\tl_range:nnn {#1} {1} {2}}~
  \mathrm{\tl_range:nnn {#1} {3} {-1}}
}

\ExplSyntaxOff

\newcommand{\arcat}[1]{\cat{Ar}(\cat{#1})}

\newcommand{\ti}[1]{\tilde{#1}}
\newcommand{\fiber}[2]{\tensor[_{#1}]{\times}{_{#2}}}

\newcommand{\cl}[1]{\mathcal{#1}}

\newcommand{\rra}{\rightrightarrows}

\newcommand{\define}[1]{\emph{#1}}

\DeclareMathOperator{\pr}{pr}

\DeclareMathOperator{\Hol}{Hol}

\newcommand{\R}{\mathbb{R}}

\newcommand{\Q}{\mathbb{Q}}

{}

\title{A singular Serre-Swan theorem via tepui fibrations}

\author{Alfonso Garmendia
	\footnote{CY Cergy Paris University, France.
		Email: 
		\texttt{garmendia-gonzalez@cyu.fr}} 
	\and 
	David Miyamoto \footnote{Queen's University, Kingston, Canada.
		Email:
		\texttt{d.miyamoto@queensu.ca}}
	\and 
	Leonid Ryvkin \footnote{Université Claude Bernard Lyon 1, France.
		Email: 
		\texttt{ryvkin@math.univ-lyon1.fr}} 
}

\begin{document}

\date{\today}

\maketitle

\begin{abstract}
The classical Serre-Swan theorem asserts that any finitely generated projective module over the algebra $C^\infty(M)$ of smooth functions of a manifold $M$ can be realized as the sections of a vector bundle over $M$. In this article, we extend this theorem beyond the projective case by introducing a notion of singular vector bundle whose sections can realize all finitely generated $C^\infty(M)$-modules, up to invisible elements. We introduce tepui fibrations as the underlying geometric objects of these singular vector bundles, and show how these tepui fibrations can model singular foliations, their holonomy groupoids, and singular subalgebroids.
\end{abstract}

\tableofcontents

\section*{Introduction}
\addcontentsline{toc}{section}{Introduction}

Given a regular foliation $\cl{F}$ of a manifold $M$, there are various approaches to understanding the geometry of the leaf space $M/\cl{F}$. Many of these are indirect, since this quotient may not be a manifold. For instance, we may consider the normal bundle $\nu \cl{F} \coloneqq TM/\cl{F} \to M$, which supports transverse constructions such as the Bott connection. Or, appealing to higher differential geometry, we may use the holonomy groupoid $\Hol(\cl{F}) \rra M$, whose arrow space is a finite-dimensional manifold that arises as a quotient of the infinite-dimensional space of foliated paths. However, when generalizing to the case of singular foliations, we run into immediate problems: the ``normal bundle'' $\nu\cl{F} = TM/\cl{F}$ no longer has constant rank, and the arrow space of the holonomy groupoid $\Hol(\cl{F})$ (as defined by Androulidakis and Skandalis \cite{AS09}) frequently fails to be a manifold. 

In this article, we demonstrate that by slightly expanding the category of manifolds to include what we call \define{tepui fibrations}, we can systematically deal with such singular objects as $\nu \cl{F}$ or $\Hol(\cl{F})$. Moreover, these singular objects retain a meaningful geometric character. In this first paper on tepui fibrations, we introduce the definition and focus on the fiberwise linear theory of tepui vector bundles and Lie tepui algebroids. A treatment of general (non-linear) tepui fibrations, including the aspects of Lie theory that will link a singular foliation with its holonomy groupoid, will appear in a subsequent article \cite{tepuiholonomy}.

Roughly speaking, a tepui fibration is a fibration $X \to M$, whose base space $M$ is a smooth manifold, and whose fibers are also smooth manifolds. The dimension of the fibers may vary, but we require that a coordinate system for one fiber can always be extended to a parametrization of the nearby fibers. We make this intuitive definition rigorous by realizing tepui fibrations as fibrations of diffeological spaces; the category $\cat{Dflg}$ of diffeological spaces, introduced by Souriau, enlarges the category of smooth manifolds, and includes all quotients and subsets of, and function spaces between, manifolds. Our tepui fibrations, formally introduced in Definition \ref{defn:tepui}, constitute a subcategory of (the arrow category of) $\cat{Dflg}$. The archetypical examples of tepui fibrations are diffeological vector bundles of the form $V/D \to M$, where $V \to M$ is a smooth vector bundle, and $D$ is a smooth \emph{singular subbundle} of $V$ (see Lemma \ref{lem:smoothsingvbtotepuis}). 

To illustrate what is gained by axiomatizing tepui fibrations, we give a generalization of the Serre-Swan theorem.\footnote{{Swan's original proof \cite{Swan62}} is in the category of topological spaces. For a proof in the smooth case, see e.g.\ \cite[Chapter 12]{nestruev}} In the classical case, recall: for a fixed manifold $M$, the global section functor
\begin{equation*}
    \Gamma\colon \{\textrm{vector bundles over }M\} \to \{\textrm{finitely generated and projective }C^\infty(M)\textrm{-modules}\}
    \end{equation*}  
    is an equivalence of categories. This means that every such $C^\infty(M)$-module is isomorphic to the space of sections of some vector bundle ($\Gamma$ is essentially surjective), and every morphism between these modules is induced by a unique morphism of vector bundles ($\Gamma$ is full and faithful). The vector bundle objects that correspond to more general $C^\infty(M)$-modules are precisely those tepui fibrations that are also diffeological vector bundles, which we call \define{VB-tepui} (vector bundle tepui).

  \begin{thmno}[Theorem \ref{thm:gensingserreswan}]
    For a fixed manifold $M$, the global section functor
    \begin{equation*}
    \Gamma\colon \{\textrm{VB-tepui over }M\} \to \{\textrm{finite type, global, and fiber-determined }C^\infty(M)\textrm{-modules}\}.
  \end{equation*}
  is an equivalence of categories.
  \end{thmno}
  We furthermore have
  \begin{itemize}
  \item The positive result that for a VB-tepui $E$, the space of sections $\Gamma(E)$ is naturally a Fr\'{e}chet space, so that $\Gamma$ takes values in the category of Fr\'{e}chet $C^\infty(M)$-modules. Our approach is to recognize that the canonical functional diffeology on $\Gamma(E)$ gives the Fr\'{e}chet structure. 
  \item The interesting negative result that, unlike in the classical case, the global section functor $\Gamma$ is no longer monoidal. This leads to Example \ref{ex:not-monoidal}, which gives two fiber-determined $C^\infty(M)$-modules whose (completed) tensor product is not fiber-determined. 
  \end{itemize}
    Singular foliations $\cl{F}$, according to their contemporary definition (cf.\ \cite{LaurentGengouxLouisRyvkin}), are instances of finite type, global, and fiber-determined $C^\infty(M)$-modules. Thus, a direct corollary of our result is a natural Fr\'{e}chet structure on $\cl{F}$, even if $\cl{F}$ is not a closed subspace of the Fr\'{e}chet space of vector fields $\mathfrak{X}(M)$.

  Continuing with singular foliations, we observe that the Lie bracket of vector fields in $\cl{F}$ equips the VB-tepui corresponding to $\cl{F}$ with the structure of a \define{Lie tepui algebroid}. Lie tepui algebroids directly generalize Lie algebroids, and arise from singular foliations, singular algebroids, and almost-Lie algebroids, and underlie $L_\infty$-algebroids. This sets up our future work \cite{tepuiholonomy}, where we will differentiate the holonomy groupoid of a singular foliation to the associated Lie tepui algebroid. 

  \subsection*{Overview}
  
In Section \ref{sec:fibrations} we recall the necessary diffeological background and introduce general and VB-tepui fibrations. We also show that VB-tepui fibrations are locally modelled by quotients of smooth vector bundles by smooth singular subbundles.

  In Section \ref{sec:singular-serre-swan} we prove Theorem \ref{thm:gensingserreswan}, the aforementioned singular Serre-Swan theorem,

  and show that the sections of a VB-tepui are naturally a nuclear Fréchet space. We close the section by showing that while the global section functor is not monoidal, we can still describe $\Gamma(E \otimes E')$ in terms of $\Gamma(E)$ and $\Gamma(E')$.

  In Section \ref{sec:algarbroids} we turn our attention to  Lie tepui algebroids. We show that singular foliations, singular subalgebroids, Lie algebroids, almost Lie algebroids, and $L_\infty$-algebroids all give rise to tepui algebroids. Unlike all these examples, however, a general tepui algebroid is not necessarily smooth along its leaves.

\subsection*{Acknowledgements}

The authors would like to thank Iakovos Androulidakis, Christian Blohmann, Lory Aintablian, Camille Laurent-Gengoux, Ruben Louis, and Marco Zambon for various discussions around the topic of this article. Special thanks go to Joel Villatoro for communicating to us his independent ongoing work on the connection between sheaves of modules and diffeological bundles.

A.\ G.\ and D.\ M.\ would both like to thank the Max Planck Institute for Mathematics in Bonn, where they completed the majority of their contribution to this paper, for providing a stimulating and collegial environment.

L.\ R.\ was supported by the DFG grant Higher Lie Theory - Project number 539126009. 

For the purpose of Open Access, a CC-BY-NC-SA public copyright license has been applied by the authors to the present document and will be applied to all subsequent versions up to the Author Accepted Manuscript arising from this submission.

\section{Tepui fibrations as diffeological spaces}
\label{sec:fibrations}
\subsection{Preliminaries on diffeological spaces}

The standard reference for diffeology is the textbook \cite{Igl13}. Let us start by recalling diffeological spaces and their natural topologies.

\begin{defn} A diffeology on a set $X$ is an assignment $\CDD$ that gives, for any natural number $k\in \KN$ and open set $U\subset \KR^k$, a subset of maps $\CDD(U)\subseteq {\rm Map}(U,X)$ called \define{plots} such that:
	\begin{itemize}
		\item (concreteness) every constant map is in $\CDD$;
		\item (pre-sheaf) for any $p\in \CDD(U)$, $k'\in \KN$, $U' \subseteq \KR^{k'}$ open, and $f\in \gi(U',U)$ the composite $p \circ f \in \CDD(U')$;
		\item (sheaf) $p\in \CDD(U)$ if and only if for any open cover $U=\bigcup_{i\in I} U_i$, we have $p|_{U_i}\in \CDD(U_i)$ for all $i \in I$.
	\end{itemize}
	A \define{diffeological space} is a pair $(X,\CDD)$ where $X$ is a set and $\CDD$ is a diffeology on $X$.

	Given diffeological spaces $(X,\CDD)$ and $(Y,\CDD')$, a function $f \colon X\fto Y$ is called \define{smooth} if $f\circ \CDD\subseteq \CDD'$. The set of smooth functions from $X$ to $Y$ will be denoted $\gi(X, Y )$. Smooth functions admitting a smooth inverse will be called \define{diffeomorphisms}.
\end{defn}

Diffeological spaces, and smooth maps between them, form the category $\cat{Dflg}$. The category $\cat{Dflg}$ is a concrete category\footnote{A concrete category over some category $\cat{A}$ is a category $\cat{C}$ equipped with a faithful functor $\cat{C} \to \cat{A}$} over $\cat{Set}$. Given a diffeological space $X$, we denote its underlying set by $|X|$.

\begin{ex}
	\label{ex:coarse-discrete-diffeologies}
	Any set $X$ has the following two diffeologies: the \define{coarse} diffeology $\cl{D}_\bullet$, where $\cl{D}_\bullet(U) \coloneqq \operatorname{Map}(U,X)$, and the \define{discrete} diffeology $\cl{D}_\circ$, where $\cl{D}_\circ(U)$ consists of only the locally constant maps into $X$. We denote $X_\bullet \coloneqq (X, \cl{D}_\bullet)$ and $X_\circ \coloneqq (X, \cl{D}_\circ)$. The coarse diffeology functor $X\mapsto X_\bullet$ is a right adjoint to the forgetful functor $|\cdot|$, while the discrete diffeology functor $X\mapsto X_\circ$ is a left adjoint to $|\cdot|$.
\end{ex}

The category $\cat{Dflg}$ is also a concrete category over topological spaces, by means of the D-topology.

\begin{defn}
	The \define{D-topology} on a diffeological space $(X, \CDD)$ is the final topology with respect to the plots of $X$. In other words, $A \subseteq X$ is open if and only if $p^{-1}(A)$ is open for all plots $p \in \CDD$.
\end{defn}

If a map $f\colon X \to Y$ between diffeological spaces is smooth, then it is continuous with respect to the D-topologies. Thus, passing to the D-topology defines a faithful functor $D\colon \cat{Dflg} \to \cat{Top}$. This functor is left adjoint to the \emph{continuous diffeology} functor, which assigns to a topological space $T$ the diffeology $\cl{D}(U) \coloneqq C(U,T)$. The D-topology functor does not have a left adjoint; for instance, in Remark \ref{rem:limits-not-preserved}, we see that $D$ does not preserve limits.

We denote by $\cat{Mfld}$ the category of smooth manifolds (Hausdorff and second-countable), with smooth maps between them. There is a functor $y\colon \cat{Mfld} \to \cat{Dflg}$ of concrete categories over $\cat{Set}$, which assigns to a manifold $M$ the diffeology $\CDD_M$ whose plots $\CDD_M(U)$ are precisely the usual smooth maps from $U$ to $M$. The functor $y$ is full and faithful, and also preserves the concrete structure over $\cat{Top}$ (i.e., the D-topology of $y(M)$ is the original manifold topology). Succinctly, $y$ is an embedding of $\cat{Mfld}$ into $\cat{Dflg}$.

The category $\cat{Dflg}$ satisfies several convenient properties: it has all small limits, all small colimits, is locally Cartesian closed, and has a weak subobject classifier. We will describe only the constructions that we need, but stress that they are all induced by the properties listed above.

\begin{defn}
	Given a diffeological space $X$ and a subset $A$, the \define{subspace diffeology} on $A$ consists of those maps $p\colon U \to A$ for which $p\colon U \to X$ is a plot.

	A map $i\colon Y \to X$ is an \define{induction} if it is injective, and $i\colon Y \to i(Y)$ is a diffeomorphism, where $i(Y)$ carries the subspace diffeology.
\end{defn}

\begin{rk}
  \label{rem:limits-not-preserved}
	The D-topology on $A$ is finer (it has more open subsets) than the subspace topology that $A$ inherits from $D(X)$, and is often strictly finer. For example, if $X = \R$ and $A = \Q$, then the plots of $\Q$ are all locally constant. Therefore, $D(\Q)$ is the discrete topology; every subset is open. On the other hand, the D-topology of $\R$ is its usual topology, and the subspace topology on $\Q$ is generated by the rational intervals; for instance, singletons are not open.
\end{rk}

The subspace topology, and D-topology of the subspace diffeology, coincide for open subsets.

\begin{lem}
	If $A$ is a D-open subset of a diffeological space $X$, then the D-topology on $A$ is exactly the subspace topology that $A$ inherits from $D(X)$.
\end{lem}
Since a diffeology naturally induces a topology, and the restriction of a diffeology to an open subset again yields a diffeology, we can talk about germs of diffeological spaces.

\begin{defn} Let $X,Y$ be a diffeological spaces and $x\in X$, $y\in Y$. We write $(X,x)$ for the \define{germ} of $X$ at $x$. A smooth map $(X,x)\to (Y,y)$ is given by an open neighbourhood $A\subset X$ of $x$ and a smooth map $A\to Y$ mapping $x$ to $y$. Two such maps are considered equivalent if there is an open neighbourhood of $x$ on which they coincide.
\end{defn}

This notion of germ allows for a concise definition of the quotient diffeology.

\begin{defn}
  \label{def:quotientdflg}
	Given a diffeological space $X$ and an equivalence relation $\sim$ on $X$, let $\pi\colon X \to X/{\sim}$ denote the quotient map. The \define{quotient diffeology} on $X/{\sim}$ consists of those maps $p\colon U \to X/{\sim}$ such that, for every $u \in U$, there is a plot $q\colon (U,u) \to X$ such that $\pi q$ and $p$ have the same germ at $u$.

	A map $s\colon X \to Y$ is a \define{subduction} if it is surjective and the natural map $Y \to X/s$ is a diffeomorphism, where the latter carries the quotient diffeology and $X/s$ denotes the equivalence classes induced by the fibers of $s$.
\end{defn}

Equivalently, the plots of $X/{\sim}$ in the quotient diffeology are those maps which locally lift along $\pi$.

\begin{rk}
	Unlike the case of the subspace diffeology, the D-topology of $X/{\sim}$ is precisely equal to the quotient topology on $D(X)/{\sim}$.
\end{rk}

\begin{rk}\label{rk:subductionsstrongepi}
	We note that inductions are precisely the strong monomorphisms, and subductions are precisely the strong epimorphisms, in $\cat{Dflg}$. As a consequence, pullbacks along subductions are subductions (cf. \cite[Proposition 3.6, Terminology 3.7]{Bloh24}).
\end{rk}

Given a diffeological space $X$ with relation $\sim$, and a subset $A \subseteq X$, there are two reasonable diffeologies on the set $A/{\sim}$, namely the quotient diffeology of the subdiffeology on $A$, or the subdiffeology of the quotient diffeology on $X/{\sim}$. The following lemma establishes a criterion for when they coincide.

\begin{lem}
	\label{lem:subquotient}

	In the following commutative diagram
	\begin{equation*}
		\begin{tikzcd}
			A \ar[r, "\iota", hook] \ar[d, "f"] & X \ar[d, "\pi", two heads] \\
			B \ar[r, "j"] & Y,
		\end{tikzcd}
	\end{equation*}
	suppose that $\iota$ is an induction, and $\pi$ is a subduction, and $\pi^{-1}(j(B)) \subseteq \iota(A)$.
	\begin{enumerate}[(a)]
		\item If $f$ is surjective, then it is a subduction.
		\item If $j$ is injective, then it is an induction.
	\end{enumerate}
\end{lem}
\begin{proof}
	\begin{enumerate}[(a)]
		\item Suppose $f$ is surjective. Let $p\colon U \to B$ be a plot. Then $jp\colon U \to Y$ is a plot, which locally lifts along $\pi$ to a plot $\widetilde{jp}$ of $X$. This plot has image in $\iota(A)$, hence $\iota^{-1}\widetilde{jp}$ is the desired local lift of $p$.
		\item Suppose $j$ is injective. Let $p\colon U \to Y$ be a plot with image in $j(B)$. Then $p$ locally lifts along $\pi$ to a plot $\ti{p}$, whose image is in $\iota(A)$. Thus, $j^{-1} p$ is locally $f \iota^{-1} \ti{p}$, which is smooth. Therefore, $j^{-1}p$ is a plot, and $j$ is an induction.
	\end{enumerate}
\end{proof}
We can now return to the case where $Y=X/{\sim}$ and $B=A/{\sim}$.  The condition $\pi^{-1}(j(B)) \subseteq \iota(A)$ then translates to requiring that elements of $A$ are equivalent only to other elements of $A$. The map $j$ can be shown to be smooth and hence, by Lemma \ref{lem:subquotient} an induction, i.e., the quotient diffeology on $A/{\sim}$ coincides with the subdiffeology on its image in $X/{\sim}$.

While in the context of manifolds, fibrations are given by surjective submersions, the picture is much less clear in the context of diffeology. 
In fact, in addition to subductions, which we already encountered in Definition \ref{def:quotientdflg}, we will need \emph{plotwise submersions}, and their surjective version \emph{local subductions}. 

\begin{defn}Let $\varphi\colon X\to Y$ be a smooth map between diffeological spaces.
	\begin{itemize}
		\item We call $\varphi$ a \define{plotwise submersion} if for any $x\in X$, and any plot $p\colon (U,0)\to (Y,\varphi(x))$, there exists a plot $q\colon (U,0)\to (X,x)$ with $\varphi q=p$ as germs at $0$.
		\item We call $\varphi$ a \define{local subduction} if it is a surjective plotwise submersion.
	\end{itemize}
\end{defn}
The term ``plotwise submersion'' is from \cite{ahmadiSubmersionsImmersionsEtale2024}. By definition, subductions are surjective maps and local subductions are subductions. Nevertheless, not all subductions are plotwise submersions, and neither the other way around.

\begin{ex}
	Any map $X \to Y$ with a global section is a subduction. Thus, for instance, the map $\R^2 \to \R$ given by $(x,y) \mapsto xy$ is a subduction. However, it is not a plotwise submersion, since the identity plot $1_{\R}\colon (\R,0) \to (\R,0)$ does not lift to a plot through $(0,0)$. Conversely, any non-surjective submersion between manifolds is a plotwise submersion that is not a subduction.
\end{ex}

Plotwise submersions are open maps for the D-topology \cite[Proposition 3.14]{ahmadiSubmersionsImmersionsEtale2024}. On the full subcategory $\cat{Mfld}$, the plotwise submersions are exactly the submersions (\cite[Corollary 3.24]{ahmadiSubmersionsImmersionsEtale2024}), and thus local subductions are the surjective submersions.

\begin{lem}
  \label{lem:pullback-submersions}
Plotwise submersions are stable under pullback.  
\end{lem}
In Remark \ref{rk:subductionsstrongepi}, we saw that subductions, being the strong epimorphisms in $\cat{Dflg}$, are stable under pullback. Plotwise submersions may not have such a convenient categorical description, so we proceed directly.
\begin{proof}
  Let $\varphi\colon X \to Y$ be a plotwise submersion, and $f\colon Z \to Y$ be a map. Let us recall that the pullback of $\varphi$ along $f$ is the projection $f^*X \to Z$, where $f^*X \coloneqq \{(z,x) \mid f(z) = \varphi(x)\}$. We claim $f^*X \to Z$ is a plotwise submersion, so let $p\colon (U,u) \to (Z,\pr(z,x))$ be a plot. Then lift $fp\colon (U,u) \to (Y,f(z))$ along $\varphi$ to a plot $q\colon (U,u) \to (X,x)$. The desired lift of $p$ is $(p,q)\colon (U,u) \to (f^*X,(z,x))$. 
\end{proof}

We also need the product and coproduct diffeologies.
\begin{defn}
	Given a family of diffeological spaces $X_\alpha$:
	\begin{itemize}
		\item the \define{product diffeology} on $\prod_\alpha X_\alpha$ consists of those maps $p = (p_\alpha) \colon U \to \prod_\alpha X_\alpha$ such that each $p_\alpha$ is a plot of $X_\alpha$.
		\item the \define{coproduct diffeology} on $\coprod_\alpha X_\alpha$ consists of those maps $p\colon U \to \coprod_\alpha X_\alpha$ such that for each $u \in U$, there is some $\alpha$ for which $p\colon (U,u) \to X_\alpha$ is a germ of a plot.
	\end{itemize}
      \end{defn}

      \subsubsection*{Locally convex spaces and their diffeology}

      We end this section by detailing the exponential objects (i.e., mapping spaces) in $\cat{Dflg}$, and compare them with some locally convex topological vector spaces. This is used in Section \ref{lem:lifting-sections}, and the tired reader may skip this passage and return to it later.
      \begin{defn}
        The \define{standard functional diffeology} on $C^\infty(X,Y)$ consists of the plots $p\colon U \to C^\infty(X,Y)$ with smooth transpose,
        \begin{equation*}
          \hat{p}\colon U \times X \to Y, \quad (u,x) \mapsto p_u(x).
        \end{equation*}
      \end{defn}
      Using the standard functional diffeology, the functor $C^\infty(X,\cdot)$ is right adjoint to the product functor $\cdot \times X$. This makes $\cat{Dflg}$ a Cartesian closed category.\footnote{Moreover, it is locally Cartesian closed, meaning that the over-categories $\cat{Dflg}\downarrow X$ are all Cartesian closed.}
      \begin{defn}\label{def:sec}
        If $\varphi\colon Y \to X$ is a smooth map, we denote by $\Gamma(\varphi)$, or $\Gamma_X(Y)$, or simply $\Gamma(Y)$, the \define{space of sections} of $\varphi$, equipped with the subdiffeology it inherits from $C^\infty(X,Y)$.
      \end{defn}
      Our key example is the space $\Gamma(V)$ of sections of a smooth vector bundle $V \to M$. We have just seen that this carries a natural diffeology. It is moreover a ``smooth $C^\infty(M)$-module'': addition and the natural multiplication by smooth functions are both smooth operations.

      On the other hand, it is more common to equip $\Gamma(V)$ with the structure of a Fr\'{e}chet vector space, and this calls for a comparison to the diffeological approach. We will assume the fundamentals of locally convex topologies, and point the reader to \cite{Trev67}, \cite{FroelKrieg88}, and \cite{gonzalessalas} as needed. The comparison begins with the following adjunctions, which we write before defining in order to set notation.
      \begin{lem}
        \label{lem:adjunctions}
        We have the following adjunctions of categories. The bottom arrows are all right adjoints to the top arrows.
        \begin{equation*}
          \begin{tikzcd}
            \mathrm{c}\cat{LCTVS} \ar[r, shift right]  & \cat{LCTVS} \ar[r, shift right, "\sigma_\infty"'] \ar[l, shift right, "c"'] & \cat{DflgVS} \ar[r, shift right] \ar[l, shift right, "\mu"'] & \cat{Dflg}. \ar[l, shift right, "F"']
          \end{tikzcd}
        \end{equation*}
      \end{lem}
      \begin{itemize}
      \item The category $\cat{LCTVS}$ is that of locally convex topological vector spaces (lctvs), and $\mathrm{c}\cat{LCTVS}$ is the subcategory of complete lctvs. The functor $c$ is the completion functor, which sends $E$ to the completion of its separation $E/\overline{\{0\}}$. The fact that $c$ is a left adjoint to the inclusion is classical, see e.g.\ \cite[Chapter 5]{Trev67} for details.
      \item The category $\cat{DflgVS}$ is the category of diffeological vector spaces. The diffeological vector space $\sigma_\infty(E)$ has for plots those maps $p\colon U \to E$ such that $l \circ p$ is smooth for all continuous linear functionals $l$. The lctvs $\mu(V)$ is equipped with the finest locally convex topology such that the continuous linear functionals on $\mu(V)$ coincide with the smooth linear functionals on $V$. This is called the \define{Mackey topology}. For the fact that $\mu$ is a left adjoint to $\sigma_\infty$, see \cite[Propositions 2.1.9 and 2.4.1]{KriegMic97}.
      \item The functor $F$ sends $X$ to the \define{free diffeological vector space} $F(X)$. Its underlying vector space is the free vector space on the set $X$, and its diffeology is generated by all maps of the form
        \begin{equation*}
          (\R \times U_1) \times \dots \times (\R \times U_n) \to F(X), \quad ((r_i,u_i))_i \mapsto \sum r_i [p_i(u_i)],
        \end{equation*}
        where $p_i\colon U_i \to X$ are plots. For the definition of $F(X)$ and the fact that $F$ is a left adjoint to the inclusion, see \cite[Section 3]{Wu15}. 
      \end{itemize}

      In particular, for the category $\cat{Frch}$ of Fr\'{e}chet spaces, the functor $\sigma_\infty$ induces the functor
      \begin{equation*}
\sigma_\infty^{\mathrm{Fr}} \colon \cat{Frch} \to \cat{Dflg}.
      \end{equation*}
      \begin{prop}
        \label{prop:frechet-vs-diffeology}
        Let $V \to M$ be a vector bundle. Then $\sigma_\infty^{\mathrm{Fr}}(\Gamma(V))$ coincides with the standard functional diffeology on $\Gamma(V)$.
      \end{prop}
  
      \begin{proof}
        To avoid bulky notation, for this proof we denote by $\Gamma_{\mathrm{Fr}} \coloneqq \sigma_\infty^{\mathrm{Fr}}(\Gamma(V))$, and by $\Gamma_{\mathrm{std}}$ the space of sections with its standard functional diffeology. First, by \cite[Lemma 30.8]{KriegMic97}, both $\Gamma_{\mathrm{Fr}}$ and $\Gamma_{\mathrm{std}}$ have the same smooth curves. We claim this implies that they have the same plots. To do this, we twice invoke Boman's theorem \cite{Bom67}, which states that a map between smooth manifolds is smooth if and only if it is smooth when pre-composed with all smooth curves into its domain.

        First, we have the equivalences
        \begin{align*}
          p\colon U \to \Gamma_{\mathrm{Fr}} \text{ is a plot } &\iff l \circ p\colon U \to \R &&\text{ is smooth for all continuous functionals }l \\
                                                                &\iff l \circ p \circ c\colon \R \to \R &&\text{ is smooth for all }l \text{ and all smooth curves }c \\
          &\iff p \circ c\colon \R \to \Gamma_{\mathrm{Fr}} &&\text{ is a plot for all }c.
        \end{align*}
Second, the map $\hat{p}\colon U \times M \to V$ is smooth (i.e., $p$ is a plot of $\Gamma_{\mathrm{std}}$) if and only if $\hat{p} \circ (c, \lambda)$ is smooth for all smooth curves $(c, \lambda)\colon \R \to U \times M$, and this composition factors as
        \begin{equation*}
          \begin{tikzcd}
            \R \ar[r, "{(1,\lambda)}"] & \R \times M \ar[r, "\widehat{p \circ c}"] & V.
          \end{tikzcd}
        \end{equation*}
        Thus, $p\colon U \to \Gamma_{\mathrm{std}}$ is a plot if and only if $p \circ c \colon \R \to \Gamma_{\mathrm{std}}$ is a plot for all $c$. We conclude that $\Gamma_{\mathrm{Fr}}$ and $\Gamma_{\mathrm{std}}$ have the same plots, hence are diffeomorphic.
      \end{proof}

      \begin{rk}
        \label{rk:d-topology-is-mackey}
        For a general lctvs $E$, there are several \emph{a priori} different notions of a smooth map $U \to E$, for instance: Michel-Bastiani (MB) smoothness, $c^\infty$-(convenient)-smoothness, and scalar-wise smoothness. We use the latter. When $E$ is a Fr\'{e}chet space, all these notions coincide, and furthermore the Mackey topology of $\sigma_\infty(E)$ is its original topology, and this is also the D-topology. For an account of these facts, see \cite[Chapter 1]{KriegMic97}.
      \end{rk}

      \begin{rk}
        One may also use the notions mentioned in the previous remark, such as MB smoothness, to define MB-smooth maps between open subsets of Fr\'{e}chet spaces. Losik \cite{Los92} showed that the natural functor from the category of open subsets of Fr\'{e}chet spaces, and MB-smooth maps between them, into $\cat{Dflg}$, is full and faithful.
      \end{rk}

\subsection{Tepui fibrations and their properties}
\label{sec:tepui-fibrations}
In this section we introduce a special class of diffeological spaces. This definition is inspired by the holonomy groupoid of a singular foliation as introduced in \cite{AS09}. There are two facts that make the arrow space of this groupoid behave almost like a manifold. One is that it is the quotient of a disjoint union of open subsets of finite-dimensional Euclidean spaces. The second is that it fibered over a manifold $M$, and for any point $p\in M$ the fiber is locally Euclidean; more precisely, there is an open subset of an Euclidean space where the quotient map from the first fact is injective on a given fiber \cite[Proposition 1.1]{debordLongitudinalSmoothnessHolonomy2013}.

\begin{defn}\label{defn:tepui} We say that a local subduction $\varphi\colon X\to M$ from a second-countable diffeological space $X$ to a smooth manifold $M$ is a \define{tepui fibration}\footnote{we explain this name in Remark \ref{rk:tepui.name}} if for any $x\in X$, there exists a plot $p\colon U\to X$ such that:
	\begin{itemize}
		\item[(Tp1)] the image of $p$ contains $x$, and $p$ is a plotwise submersion ($p$ is a \define{chart} at $x$);
		\item[(Tp2)] the restriction
		      \begin{equation*}
			      p\colon U_{\varphi(x)} \coloneqq (\varphi p)^{-1}(\varphi(x)) \to X_{\varphi(x)}
		      \end{equation*}
		      is injective ($p$ is \define{minimal} at $x$);
		\item[(Tp3)] the fibers $X_m = \varphi^{-1}(m)$ are Hausdorff and second-countable for all $m\in M$.
	\end{itemize}

	For the rest of the paper, plots satisfying the first condition above are called \define{charts} of $X$, and those satisfying the first two conditions are called \define{minimal charts}.
\end{defn}

The following property of $\varphi$-fibers in a tepui fibration motivates our definition.
\begin{prop}
	\label{prop:fibers-are-manifolds}
	Let $\varphi\colon X \to M$ be a tepui fibration. Then the $\varphi$-fibers are manifolds. More precisely, for a minimal chart $p$ at $x$, writing $m \coloneqq \varphi(x)$, the fiber $U_m$ is a submanifold of $U$, and the restriction $p\colon U_m \to X_m$ is a diffeomorphism onto a D-open subset of $X_m$.
\end{prop}
\begin{proof}
	Fix a minimal chart $p\colon U \to X$ at $x$. The composite $\varphi p\colon U \to M$ is a plotwise submersion between manifolds, hence is a submersion. Therefore, $U_m$ is an embedded submanifold of $U$ (and its sub-diffeology coincides with its manifold structure).

	Let $q\colon W \to X_m$ be a plot, and suppose $q(r) = p(u)$. Since $X_m$ carries the subspace diffeology, $q\colon W \to X$ is a plot, and since $p$ is a plotwise submersion, there is some lift $\ti{q}\colon (W,r) \to (U,u)$ of $q$ along $p$. But the image of $\ti{q}$ must be in $U_m$, since $p(\ti{q}(W)) \subseteq q(W) \subseteq X_m$. The map $\ti{q}\colon W \to U_m$ is smooth because $U_m$ carries the subspace diffeology. Thus, we have shown that $p\colon U_m \to X_m$ is a plotwise submersion.

	Then $p$ has D-open image, and the map $p\colon U_m \to p(U_m)$ is a subduction. It is injective because $p$ is minimal at $x$, hence it is a diffeomorphism. We conclude that $X_m$ is locally Euclidean. It is second-countable and Hausdorff by assumption, hence $X_m$ is a manifold.
\end{proof}

Tepui fibrations form a full subcategory of the arrow category $\arcat{Dflg}$. The import of this definition is that $X$ is a fibration over a manifold whose fibers are manifolds of perhaps varying dimension. Surjective submersions between smooth manifolds are tepui fibrations. In particular, any smooth manifold is a tepui fibration over a point.

\begin{rk}\label{rk:tepui.name}
	The term ``tepui'' was not germane to mathematics. Rather, it is a geological term: a \define{tepui} is a particular type of mesa found in South America, distinguished by its striking cliffs and surrounded by a savanna, such as in this Venezuelan example below.  Given a tepui fibration $\varphi\colon X \to M$, we believe that $X$ resembles a tepui mountain rising over the smooth savanna $M$.
	\begin{figure}[h]
		\centering
		\includegraphics[scale=0.5]{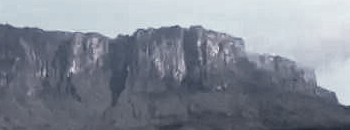}
		\caption{Roraima tepui (photo by A.G.)}
	\end{figure}

\end{rk}

\begin{rk}

	If we omit (Tp2) and (Tp3) in Definition \ref{defn:tepui}, we obtain a notion which was called \emph{generated by local subductions} in the article \cite{AndroulidakisZambonIntegrationOfSingular2020}. Since this is essentially half of our conditions for a tepui fibration, we are tempted to call such spaces \emph{hemitepui} (hemi- means half, like in hemisphere. We think of the tepui without the savanna). Note that to be hemitepui is a property of a space, whereas to be tepui fibration is a property of a map. Smooth manifolds are hemitepui. Countable disjoint unions and finite products of hemitepui remain hemitepui. Infinite-dimensional spaces (section spaces of fiber bundles, diffeomorphism groups) are not hemitepui. We will not pursue hemitepui beyond this remark.

\end{rk}

\begin{ex}
	Consider a cylinder $S^1\times \KR$ with the following relation: $(\theta,x)\sim(\sigma,y)$ if and only if $x=y\neq 0$. Then $\pr_2\colon (S^1\times \KR)/{\sim} \to \KR$ is a tepui fibration.

	On the other hand, if the relation on the cylinder is $(\theta,x)\sim_2(\sigma,y)$ if and only if $x=y= 0$, the map $\pr_2\colon(S^1\times \KR)/{\sim_2}\fto \KR$ is not an tepui fibration. For $x=0\in \KR$, there is no minimal chart, since charts, being plotwise submersions, must be open maps.
\end{ex}

For the next example, refer to Definition \ref{def:singular-foliation} for the definition of a singular foliation and to  \cite{AS09} for the definition of its holonomy groupoid. 
\begin{ex}
  \label{ex:holonomy-groupoid}
  Given a smooth manifold $M$ and singular foliation $\CF$ on $M$, the holonomy groupoid $\CH(\CF)\soutar M$, originally introduced in \cite{AS09}, satisfies:
	\begin{itemize}
		\item $\CH(\CF)\soutar M$ is a diffeological groupoid, that is, a groupoid internal to $\cat{Dflg}$. We denote its structure maps with bold font, i.e., $\bm{s}, \bm{t}$ for source and target, and $\bm{u}$ for unit.
		\item The diffeological space $\CH(\CF)$ is defined as a quotient of a disjoint union of open subsets of Euclidean spaces. Precisely, for any element $\phi\in \CH(\CF)$ there is an open set $U$ in a Euclidean space and a subduction $\varphi_U\colon (U,u)\fto (\CH(\CF),\phi)$. In the original article \cite{AS09}, the maps $\varphi_U\colon (U,u)\fto (\CH(\CF),\phi)$ are called \define{bisubmersions} of $\CF$. 
		\item In \cite[Proposition 2.10]{AS09} it is proven that for $\phi\in \CH(\CF)$ near the identity bisection $M\fto \CH(\CF)$, and bisubmersions $\varphi_{U_1}\colon (U_1,u_1)\fto (\CH(\CF),\phi)$ and $\varphi_{U_2}\colon (U_2,u_2)\fto (\CH(\CF),\phi)$ (which we know are subductions from the previous point), there is a smooth map $(U_1,u_1)\fto (U_2,u_2)$ commuting with $\varphi_{U_1}$ and $\varphi_{U_2}$. This proves that near the identity bisection, the bisubmersions are plotwise submersions, therefore they are charts.
		\item Moreover, in \cite[Proposition 1.1]{debordLongitudinalSmoothnessHolonomy2013} it is proven that for any point $m\in M$ there is a bisubmersion $\varphi_U\colon (U,u)\fto (\CH(\CF), \bm{u}(m))$ such that $\varphi_U|_{U_m}\colon U_m\fto \CH(\CF)_m=\bm{s}^{-1}(\{m\})$ is injective. Therefore, there are minimal charts for elements close enough to the identity bisection.
	\end{itemize} 
Using that $\CH(\CF)$ is a diffeological groupoid, that for any $\phi\in \CH(\CF)$ there is a bisection passing through it, and that right translations are smooth, one gets minimal charts for $\bm{s}$ at any element of $\CH(\CF)$. We conclude that $\bm{s}\colon \CH(\CF)\fto M$ is a tepui fibration. With a similar argument, or using the inverse map, we see that $\bm{t}\colon \CH(\CF)\fto M$ is also a tepui fibration.
	
\end{ex}

Here we gather some basic facts about tepui fibrations.

\begin{lem}\label{lem:generic-regularity}
  Let $\varphi\colon X \to M$ be a tepui fibration.  The \emph{minimal rank} map
  \begin{equation*}
    \mrank \colon m \mapsto \min \{\dim_x X_m \mid x \in X_m\}
  \end{equation*}
  is upper semicontinuous. Furthermore, there is an open and dense subset of $M$ on which $\mrank$ is locally constant.
  
\end{lem}
Here $\dim_x X_m$ denotes the dimension of the connected component of $X_m$ containing $x$. Also recall that $\mrank$ being upper semicontinuous is equivalent to the sets $\{m \mid \mrank(m) \leq k\}$ being open in $M$. We call a point, $m$, \define{regular} if there is an open neighbourhood of $m$ on which $\mrank$ is constant. Thus, our lemma states that the set of regular points is open and dense.

\begin{proof}
  Fix a non-negative integer $k$, and suppose that $\mrank(m_0) \leq k$. Say this is witnessed by $\dim_{x_0} X_{m_0} = \mrank(m_0)$, and fix a minimal chart $p\colon U \to X$ at $x_0$. Without loss of generality, we may assume that $\varphi p$ has connected fibers. At any other point $m$ in the open subset $\varphi p(U) \subseteq M$, the restricted map $p\colon U_m \to X_m$ is a plotwise submersion between manifolds (examine the proof of Proposition \ref{prop:fibers-are-manifolds}), hence a submersion. Therefore, at any $x$ in the open subset $p(U_m) \subseteq X_m$, we have
  \begin{equation*}
    \mrank(m) \leq \dim_x X_m \leq \dim U_m = \dim_{x_0} X_{m_0} = \mrank(m_0) \leq k.
  \end{equation*}
  For the second assertion, observe that since $\mrank$ takes values in $\mathbb{Z}_{\geq 0}$, on any given open subset $O$ of $M$ it attains a minimum value, say $k$. Then any point $m$ in $O$ with $\mrank(m) = k$ is regular, of which one must exist.
\end{proof}

\begin{lem}\label{lem:workdonkeylemma}
	Let $p\colon (U,u)\to (X,x)$ be a minimal chart for a tepui fibration $\varphi\colon X \to M$. Assume $\phi\colon (U,u)\to (U,u)$ lies over $1_X$, so that the following diagram commutes:
	\begin{equation*}
		\begin{tikzcd}
			(U,u) \ar[r, "\phi"] \ar[d, "p"] & (U,u) \ar[d, "p"] \\
			(X,x) \ar[r, "1_X"] & (X,x).
		\end{tikzcd}
	\end{equation*}
	Then $T_u\phi$ is invertible, thus $\phi$ is a diffeomorphism near $u$.
\end{lem}
\begin{proof} The composition $\varphi p\colon U \to M$ is a submersion, so we may realize $(U,u)\cong (M\times W, (m,w))$ with $m=\varphi(x)$. We observe that $T_uU= T_wW \oplus T_{m}M$. Since $\phi$ is lies over the identity on $M$ we know that
        \begin{equation*}
          T_u\phi=
          \begin{pmatrix}
            1_{T_{m}M} & B            \\
            0 & A
          \end{pmatrix}
        \end{equation*}
	where $A\colon T_wW\to  T_wW$, and $B\colon T_{m}M\to T_wW$. The map $p\colon U_m \cong \{m\} \times W \to X_m$ is injective, hence  $\phi|_{\{m\}\times W} = 1_W$. This implies that  $A=1_{T_wW}$, which gives the invertability of $T_u\phi$ and hence the local invertability of $\phi$.
\end{proof}

We can use this to deduce how any chart factors through a minimal one:

\begin{cor}
  \label{cor:localequiv}
  Let $q\colon (W,w)\to (X,x)$ be a chart and $p\colon (U,u)\to (X,x)$ a minimal chart. Then any lift $\tilde q$ of $q$ along $p$ is a submersion near $w$. Similarly, any lift $\tilde p$ of $p$ along $q$ is an immersion near $u$.
\end{cor}

\begin{proof}  Consider to $\tilde q\circ \tilde p\colon (U,u)\to (U,u)$. It lies over $1_X$, since $p\tilde q\tilde p=q\tilde p=p$ by construction. Hence, we can apply Lemma \ref{lem:workdonkeylemma} and obtain that $\tilde q\tilde p$ is invertible near $u$, whence $\tilde q$ is a submersion and $\tilde p$ an immersion.

\end{proof}

Pulling back tepui fibrations along smooth maps from manifolds yields new tepui fibrations.

\begin{lem}
  \label{lem:pullback-of-tepui}
  Let $\varphi\colon X \to M$ be a tepui fibration, and let $f\colon N \to M$ be a smooth map between manifolds. Then $f^*X \to N$ is a tepui fibration.
\end{lem}
\begin{proof}
  Fix a chart $\varphi\colon U \to X$, and assume $\varphi$ is minimal at $m$, so that $\varphi \colon U_m \to X_m$ is injective. By the pasting lemma, both small rectangles and the large rectangle in the diagram below are pullbacks:
  \begin{equation*}
    \begin{tikzcd}
      f^*U \ar[r, "\ti{p}"] \ar[d] & f^*X \ar[r, "\ti{\varphi}"] \ar[d] & N \ar[d, "f"] \\
      U \ar[r, "p"] & X \ar[r, "\varphi"] & M.
    \end{tikzcd}
  \end{equation*}
  By Lemma \ref{lem:pullback-submersions}, the map $\ti{\varphi}$ is a local subduction, and the map $\ti{p}$ is a plotwise submersion. Furthermore, for $n \in N$ with $f(n) = m$, the map $\ti{p}$ is
  \begin{equation*}
    \ti{p}\colon \{n\} \times U_m \to \{n\} \times X_m, \quad (n,u) \mapsto (n, p(u)),
  \end{equation*}
  which is injective because $p\colon U_m \to X_m$ is injective. Thus, $\ti{p}$ is a minimal chart at $n$. We also see that the fibers of $\ti{\varphi}$ identify with the fibers of $\varphi$, thus are Hausdorff and second-countable. We conclude that $\ti{\varphi}\colon f^*X \to N$ is indeed a tepui fibration.
\end{proof}

\subsection{Tepui vector bundles}

\label{sec:vb-tepui}

In this subsection, we will describe the ``vector bundle'' objects in our context. We start by recalling:
\begin{defn} A \define{(diffeological) bundle of vector spaces} is given  by diffeological spaces $E$ and $X$ and maps:
  \begin{alignat*}{2}
    &\text{\define{projection}} & p&\colon E \to X \\
    &\text{\define{addition}} & +&\colon E \fiber{}{X} E \to E \\
    &\text{\define{zero section}} & 0&\colon X \to E \\
    &\text{\define{scalar multiplication}}\quad  & \cdot &\colon \R \times E \to E
  \end{alignat*}
  such that $0$, $+$, and $\cdot$ are compatible with the projection $p$ and satisfy all the axioms of a vector space on the $p$-fibers.  
\end{defn}

Crucially, we do not assume that bundles of vector spaces admit local trivializations.

\begin{rk}
  We adopted the naming from \cite[Terminology 2.2]{Bloh24} and note that bundles of vector spaces also appear under the names vector pseudo-bundles \cite{pervova}, diffeological vector spaces over a base \cite{chriswu16}, and regular vector bundles \cite{vincent2008diffeological}.
\end{rk}

\begin{ex}
  A vector bundle $V \to M$ in $\cat{Mfld}$ is a bundle of vector spaces. Conversely, if $V \to M$ is a bundle of vector spaces, and $V$ and $M$ are manifolds, one may prove that the projection $V \to M$ admits local trivializations, in which case we have a vector bundle in $\cat{Mfld}$.

  If $r\colon W \to V$ is a vector bundle morphism, then both $r(W) \to M$ and $V/r(W) \to M$ are bundles of vector spaces. They are vector bundles if and only if $r$ has constant rank.
\end{ex}

\begin{defn}
	A bundle of vector spaces $p \colon E\to M$, over a manifold $M$, is called \define{VB-tepui} if $p$ is a tepui fibration.
\end{defn}

One source of examples for tepui fibrations arises from smooth singular subbundles. Let $V\to M$ be a vector bundle in $\cat{Mfld}$. We call $D\subset V$ a \define{singular subbundle} in $V$ if $D_m\subset V_m$ is a vector subspace for any $m\in M$. We call such a subbundle \define{smooth} if for any $m\in M$ and $v\in D_m$ there is a smooth section $\sigma\colon M\to V$ with $\sigma(x)=v$ and $\sigma(M) \subseteq D$. Smooth singular subbundles are bundles of vector spaces.

\begin{ex}
  The image of any morphism $r\colon W \to V$ of vector bundles over $M$ is a smooth singular subbundle of $V$. On the other hand, the kernel of $r$ may or may not be a smooth singular subbundle.
\end{ex}

\begin{lem}\label{lem:smoothsingvbtotepuis}
	Let $\pi_V\colon V\to M$ be a vector bundle over a smooth manifold and $D\subset V$ a smooth singular subbundle. Then for $E\coloneqq V/D$ with the quotient diffeology, the projection $\pi_E\colon E \to M$ is a VB-tepui.
\end{lem}
\begin{proof}
	Since both $V\to M$ and $V \to E$ are local subductions, so is $E\to M$. Second countability and Hausdorffness of the fibers of $\pi_E$ are also automatic. To prove that $\pi_E$ is a tepui fibration, we prove that it admits minimal charts. We write $[v]$ for the equivalence class of $v\in V$ inside $E$.
        
	Let $[v]\in E$ and $m \coloneqq \pi_V(v)$. Choose a collection of sections $\cl{S} \coloneqq \{\sigma_1,\dots, \sigma_k, \alpha_1,\dots, \alpha_l\}$ of $V \to M$ with the following properties:
        \begin{itemize}
        \item $\{\sigma_1(m),\dots, \sigma_k(m)\}$ is linearly independent and projects to a basis of $E_m$ (which is a vector space since the quotient is fiberwise linear);
        \item $\alpha_i(M) \subseteq D$ for all $1\leq i \leq l$, and $\{\alpha_1(m),\dots, \alpha_l(m)\}$ is a basis of $D_m$.
        \end{itemize}
        Choose a neighbourhood $U$ of $M$ for which the collection $\cl{S}$ is a frame of $V|_U \to U$. Then we have
        \begin{equation*}
          \begin{tikzcd}
            U \times \R^n \ar[r] \ar[d,] & V_U \ar[d] \\
            U \times \R^k \ar[r, "P"] & E_U
          \end{tikzcd}
          \quad
          \begin{tikzcd}
            (u,\lambda_1,\dots,\lambda_k,\rho_1,\dots,\rho_l) \ar[r, mapsto] \ar[d, mapsto] & \sum \lambda_i \sigma_i(u) + \sum \rho_i \alpha_i(u) \ar[d] \\
            (u,\lambda_1,\dots, \lambda_k) \ar[r, mapsto, "P"] & \left[\sum \lambda_i \sigma_i(u)\right].
          \end{tikzcd}
        \end{equation*}
        The top map is a diffeomorphism, and the downward maps are local subductions. Thus, $P$ is a local subduction. Furthermore, $P$ restricted to the fiber of $\pi_E \circ P$ over $m$ is injective, so $P$ is a minimal chart around $[v]$.

        Let us now discuss the VB-structure. The zero section, scalar multiplication, and addition on $E$ are smooth because they are directly induced by those on $V$ and the following three maps are subductions:
        \begin{equation*}
          V \to E, \quad \R \times V \to \R \times E, \quad V \fiber{}{M} V \to E \fiber{}{M} E.
        \end{equation*}

\end{proof}

\begin{cor}
  \label{cor:cokernel}
  Let $r\colon W\to V$ be a vector bundle morphism over $M$. Then $r(W)$ is a smooth singular subbundle, hence $V/r(W)$ is a VB-tepui over $M$.
\end{cor}

\begin{ex} As a particular instance of the above, consider the trivial vector bundle $V = \mathbb R^2 \to \R$, and the bundle endomorphism $r(x, v) = x \cdot v$. Then the corresponding VB-tepui $E = V/r(V)$ is depicted below.
  \begin{equation*}
    \begin{tikzpicture}
      \clip(-4,-1) rectangle (4,1);
      \draw [line width=2pt] (-4,0) -- (0,0);
      \draw [line width=2pt] (0,0) -- (4,0);
      \draw [line width=1pt, color=red] (0,-1) -- (0,1);
      \draw (0,0) node [] {)(};
    \end{tikzpicture}    
  \end{equation*}

  Points in the thin red line (the fiber over $0$) are inseparable from the black rays (the image of the zero-section). The bundle projection collapses the thin red line to a point.

  One could be tempted to think of this as a cross $\{xy = 0\} \subset \mathbb R^2$ sitting over the horizontal axis; however, diffeologically it is a very different space. For instance, $E$ admits sections through any of its points, whereas the cross over the horizontal axis only carries sections through the points $(x,0)$.
\end{ex}

Another way to obtain VB-tepui is by creating new ones from old ones.

\begin{lem}\label{lem:pullbackvb} Let $E\to M$ be a VB-tepui.
\begin{itemize}
	\item Let $N$ be a manifold and $f\colon N\to M$ be a smooth map. Then $f^*E\to N$ is VB-tepui.
	\item Let $E'\to M$ be another VB-tepui. Then $E\oplus E'\to M$ is VB-tepui.
\end{itemize}
\end{lem}

\begin{proof} 
  We may deal with the first item swiftly. By Lemma \ref{lem:pullback-of-tepui}, we know that $f^*E\to N$ is a tepui fibration. Moreover, all structure maps of the bundle of diffeological vector spaces $E\to M$ pull back to the corresponding structure maps on $f^*E\to N$.

Similarly, for the second item, direct sums of bundles of diffeological vector spaces are again bundles of diffeological vector spaces, so we need only show that $E\oplus E'\to M$ is a tepui fibration. For this, we first observe that products of tepui fibrations are tepui, i.e., $E\times E'\to M\times M$ is a tepui fibration. Then $E\oplus E'=E\times_ME'$ is simply the pullback of $E\times E'$ along the diagonal $M\to M\times M$, hence is again tepui by Lemma \ref{lem:pullback-of-tepui}.
\end{proof}

It turns out that the converse of Lemma \ref{lem:smoothsingvbtotepuis} is true locally.

\begin{prop}\label{prop:tepui.is.quotient} Let $E\to M$ be a VB-tepui. Then, about each $m\in M$ there is an open neighbourhood $U$, a vector bundle $V\to U$, and a smooth singular subbundle $D$ of $V$, such that $E_U\cong V/D$.
\end{prop}

\begin{proof}
	Let $m\in M$. Let $e_1,\dots ,e_n$ be a basis of $E_m$. Since $E\to M$ is a local subduction, and $M$ is a manifold, there are local sections $\sigma_i\colon (M,m)\to (E,e_i)$. Let $U$ be the intersection of the domains of $\sigma_i$. We define:
        \begin{equation*}
          q\colon U\times \mathbb R^n\to E_U,\quad q(u,\lambda_1,...,\lambda_n) \coloneqq \sum \lambda_i\sigma_i(u).
        \end{equation*}
        We claim that $q$ is a local subduction. Since by the use of scalar multiplication we can move any point arbitrarily close to the origin of its fiber, it suffices to show that $q\colon (U \times \R^n,(m,\vec{0})) \to (E_U, 0_m)$ is a germ of a local subduction. Let $q _{\text{min}}\colon (U',0)\to (E,0_m)$ be a minimal chart. There exists a lift $\tilde q\colon (U \times \R^n, (m,\vec{0})) \to (U',0)$ of $q$ along $q_{\text{min}}$. The lift $\tilde q$ is a local diffeomorphism at $(m,\vec{0})$, because (at the level of germs)
        \begin{align*}
          &\ti{q}|_{\{m\} \times \R^n}\colon \R^n \to U'_m \text{ is a local diffeomorphism at } \vec{0}, \text{ and} \\
          &\ti{q}|_{U \times \{0\}}\colon U \to \ti{q}(U) \text{ is a local diffeomorphisms at } m.
        \end{align*}
        The first assertion uses minimality of $q_{\text{min}}$ at $0_m$ and the fact that $\{e_1,\dots, e_n\}$ is a basis of $E_m$ The latter uses the fact that $\ti{q}|_{U \times \{0\}}$ is a section of the submersion $\pi_E \circ q_{\text{min}}$. To conclude that $T_{(m,\vec{0})}\ti{q}$ is invertible, use that $T_0U' = T_0 U'_m \oplus T_0 \ti{q}(U)$. 

        Since $q = q_{\text{min}} \circ \ti{q}$ (locally) factors as the composition of two local subductions, the germ $q\colon (U \times \R^n, (m, \vec{0})) \to (E_U, 0_m)$ is a local subduction, and so is $q$. Let $D\subset V \coloneqq U \times \R^n$ be the kernel of $q$. Then $V/D \cong E_U$ is immediate, and it remains to show that $D$ is a smooth singular subbundle of $V$. Let $d \in D_u$. We lift the zero section $0\colon (U,u)\to (E_U,0_u)$ along the local subduction $q\colon (V,d)\to (E,0)$, through the point $d$. The resulting map is a section of $D$ through $d$, and hence $D$ is smooth.
\end{proof}

When the $\mrank$ of $E$ is globally bounded, the construction can be globalized, resulting in a complete description of such VB-tepui.

\begin{cor}\label{cor:globfin} Let $E\to M$ be a VB-tepui, such that the dimension of its fibers $E_m$ is globally bounded. Then there exists a vector bundle $V$ and a smooth singular subbundle $D$ such that $E=V/D$.
\end{cor}

\begin{proof} The proof uses the Brouwer-Lebesgue paving principle and can be directly carried over from the proof of \cite[Proposition 1.2.37]{LaurentGengouxLouisRyvkin}. We will sketch the proof here for completeness.
	
Let $n \coloneqq \dim(M)$ and $N \coloneqq \max_{m\in M}\mrank(m)=\max_{m\in M}\dim(E_m)$. Since the rank of $E$ is upper semicontinuous, for any $m$ there is an open neighbourhood $U_m$ and sections $e^m_1,\dots,e^m_n$ spanning $E_{U_m}$. By Ostrand's theorem (\cite[Lemma 3]{ostrandCoveringDimensionGeneral1971}), there is a cover $\{V_{i,\alpha}\}_{i\in \{1,\dots,n+1\},\alpha\in \mathbb N}$ of $M$ with the following properties:
\begin{itemize}
	\item The cover refines $\{U_m\}_{m\in M}$
	\item For fixed $i$, the  $V_{i,\alpha}$ and $V_{i,\beta}$ are disjoint.
\end{itemize}
In particular, the cover $\{W_i\coloneqq \bigcup_\alpha V_{i,\alpha}\}$ contains $n+1$ sets, each of which is the disjoint union of countably many sets (each of which is inside some $U_{m}$). Due to the disjointedness, we can combine the generators of the $U_m$ to generators $e^i_1,\dots,e^i_N$ of $W_i$. We now pick a partition of unity $\rho^i$ subordinate to $W_i$. The $\rho^ie^i_j$ extend by zero to sections $\tilde e^i_j$ defined on $M$. The collection 

\begin{equation*}
  \{ \tilde e^i_j \mid i\in \{1,\dots,n_1\}, j\in\{1,\dots,N\} \}
\end{equation*}
generates $E$. Hence, we have a map $M\times \mathbb R^{(n+1)N}\to E$. The kernel of this map forms a smooth singular subbundle, which concludes the proof.
\end{proof}

To end, we observe that every VB-tepui of globally bounded rank is of the form in Corollary \ref{cor:cokernel}

\begin{cor}
	Let $E \to M$ be a VB-tepui of globally bounded rank. Then there exists a vector bundle morphism $r\colon W \to V$ over $M$ such that $E \cong V/r(W)$.
\end{cor}

\begin{proof}
	By Corollary \ref{cor:globfin}, it suffices to show that, fixing a vector bundle $V \to M$ and a smooth singular distribution $D \subseteq V$, there is some vector bundle morphism $r\colon W \to V$ with $r(W) = D$. This, in turn, is a consequence of the following fact, proved separately as the main result in \cite{DragLeeParRic12} and \cite{Sus08}: smooth singular subbundles are globally finitely generated. In other words, there exists a finite collection of sections $\{\sigma_1,\ldots, \sigma_N\} \subseteq \Gamma(V)$ such that $D = \operatorname{span}(\sigma_1,\ldots, \sigma_N)$. Then we may take $W \coloneqq M \times \R^N$, and
	\begin{equation*}
		r\colon M \times \R^N \to V, \quad r(m,\lambda_1,\ldots, \lambda_N) \coloneqq \sum_{i=1}^N \lambda_i\sigma_i(m).
	\end{equation*}
\end{proof}

\section{A singular Serre-Swan theorem}
\label{sec:singular-serre-swan}
\subsection{A Serre-Swan theorem for finitely generated modules}

Vector bundles over a smooth manifold $M$ correspond to finitely generated projective modules over $C^\infty(M)$ by the classical Serre-Swan theorem. Tepui fibrations allow us to extend this correspondence to singular settings. We begin with the section functor $\Gamma$ from the category of diffeological bundles of vector spaces over $M$ to the category of $C^\infty(M)$-modules. Let us recall that $\Gamma$ was defined in greater generality in Definition \ref{def:sec}. On morphisms, it is given by
\begin{equation*}
  \Gamma(f\colon E \to E') \coloneqq f_*\colon\Gamma(E) \to \Gamma(E').
\end{equation*}

In order to formulate our singular Serre-Swan theorem, we need the following notion from \cite{nestruev}.
\begin{defn} Let $Q$ be a $C^\infty(M)$-module. An element $\sigma$ of $Q$ is called \define{invisible}, if 
$$
[\sigma]_m=0 \in \frac{Q}{I_mQ}
$$
for all $m\in M$, where $I_m$ denotes the ideal of smooth functions vanishing at the point $m$. We call the space of invisible elements $\mathrm{inv}(Q)$. The module $Q$ is called \define{fiber-determined} if it $\mathrm{inv}(Q)=\{0\}$, i.e., if has no non-trivial invisible elements. 
\end{defn}
The invisible elements of any module $Q$ form a submodule that can be equivalently described as $\mathrm{inv}(Q)=\bigcap_{m\in M}I_mQ$. In particular, we can form the quotient module $\frac{Q}{\mathrm{inv}(Q)}$, which turns out to always be fiber-determined. The assignment $Q\mapsto \frac{Q}{\mathrm{inv}(Q)}$ extends to a functor from the category of $C^\infty(M)$-modules to the category of fiber-determined $C^\infty(M)$-modules. We call this functor \define{fiber-determination}.

\begin{rk}\label{rem:geomvsfibdet} Fiber-determined modules and invisible elements were introduced in \cite{nestruev}. However, there, fiber-determined modules were referred to as \emph{geometric modules} and consequently \emph{fiber-determination} was called \emph{geometrization} there.
They later reappeared under the name fiber-determined modules in \cite{villa}.  We use the terminology fiber-determined here, since it is more self-explanatory. 
\end{rk}

Intuitively, to be fiber-determined means that an element of $Q$ vanishes if it vanishes at all points of $M$. Examples of fiber-determined modules include modules of sections of vector bundles, of tepui fibrations, and all submodules of these. Let us see what can go wrong in a simple example.

\begin{ex}
  \label{ex:not-geometric}
  Consider $M=\mathbb R$. We write $I_0$ for the ideal of functions vanishing at the origin. Then the $C^\infty(\R)$-module $Q\coloneqq C^\infty(\KR)/I_0^2$ is not fiber-determined: the class of the function $f(x)=x$ in this module is not zero, but its value in $Q/I_mQ$ is zero for all $m \in\KR$. To see that, we can calculate
  \begin{equation*}
	Q/I_mQ=\frac{C^\infty(\KR)/I_0^2}{I_m (C^\infty(\KR)/I_0^2)}=\frac{C^\infty(\KR)/I_0^2}{I_m C^\infty(\KR)}=C^\infty(\KR)/(I_0^2+I_m).    
  \end{equation*}
  The ideal $I_m$ is generated by $x-m$, and the ideal $I_0^2$ is generated by $x^2$. When $m\neq 0$, we have
  \begin{equation*}
    \frac{1}{m^2}x^2-\frac{1}{m^2}(x-m)(x+m)=1,
  \end{equation*}
thus $I_0^2+I_m=C^\infty(\KR)$, i.e., $Q/I_mQ=0$, and hence necessarily $[f]_m=0\in Q/I_mQ$. When $m=0$, we have $I_m+I_0^2=I_0$, and $f\in I_0$, so $[f]_0=0\in Q/I_0Q$. However, $f\not\in I_0^2$, hence the class of $f$ in $Q$ does not vanish. Consequently, $Q$ is not fiber-determined.
\end{ex}

We can now work towards a preliminary version of our singular Serre-Swan theorem.
\begin{lem}\label{lem:esssurj}Let $Q$ be a finitely generated, fiber-determined $C^\infty(M)$-module. Then there is a VB-tepui $E \to M$ of globally bounded rank with $\Gamma(E)\cong Q$.
\end{lem}
\begin{proof}
  Since $Q$ is finitely generated, there exists some vector bundle $V \to M$ and a surjection $f\colon \Gamma(V)\to Q$. Namely, if $X_1,\dots, X_N$ are generators of $Q$, we may take the trivial bundle $V \coloneqq M \times \R^N$, and $f(\sigma_1,\dots, \sigma_N) \coloneqq \sum_i \sigma_i X_i$. Let
  \begin{equation*}
  K \coloneqq \ker(f) \subseteq \Gamma(V) \quad\text{and}\quad D \coloneqq \{k(m) \mid k\in K,\ m\in M\} \subseteq V.  
  \end{equation*}
   By definition, $D$ is a smooth singular subbundle of $V$ and, $K\subset \Gamma(D)$. We show that in fact $K = \Gamma(D)$, using the fact that $Q$ is fiber-determined. Let $d \in \Gamma(D)$. For $m \in M$, fix $k \in K$ with $d(m) = k(m)$. We track where $\bar{k}$ and $k$ are sent in the diagram below: 

\begin{equation*}
\begin{tikzcd}
	{d,k} &&&& {f(d),f(k)} \\
	& {\Gamma(V)} && Q \\
	& {\Gamma(V)/I_m\Gamma(V)} && {Q/I_mQ} \\
	{d(m),k(m)} &&&& {[f(d)],0}
	\arrow[dashed, from=1-1, to=1-5]
	\arrow[dashed, from=1-1, to=4-1]
	\arrow[dashed, from=1-5, to=4-5]
	\arrow[from=2-2, to=2-4]
	\arrow[from=2-2, to=3-2]
	\arrow[from=2-4, to=3-4]
	\arrow[from=3-2, to=3-4]
	\arrow[dashed, from=4-1, to=4-5]
\end{tikzcd}
\end{equation*}

But $d(m) = k(m)$, so $[f(d)] = 0$ in $Q/I_mQ$, hence by fiber-determinedness $f(d)=0$ and $d \in K$. Since $D$ is a smooth singular subbundle of $V$, the fibration $E=V/ D\to M$ is a VB-tepui. By definition, a section of $E=V/ D$ can be locally lifted to a section of $V$, and using a partition of unity on $M$, it can even be globally lifted, thereby the quotient map $\Gamma(V)\to\Gamma(E)$ is a surjection. On the other hand, a section in $\Gamma(V)$ is mapped to zero in $\Gamma(E)$, if and only if it is in $D$ at all points, i.e., if it is a section of $\Gamma(D)$. In total, we obtain the isomorphisms of modules:
\begin{equation*}
  \Gamma(E)=\Gamma(V/D)\cong\Gamma(V)/\Gamma(D)=\Gamma(V)/K\cong Q.
\end{equation*}

\end{proof}

On the other hand, from Corollary \ref{cor:globfin}, any VB-tepui of globally bounded rank induces a finitely generated fiber-determined $C^\infty(M)$-module, via the smooth sections functor. Altogether, we obtain:

\begin{prop}
  \label{prop:bounded-serre-swan}
  The section functor $\Gamma$ induces an equivalence of categories between the category of VB-tepui over $M$ of globally bounded rank, and finitely generated fiber-determined $C^\infty(M)$-modules.
\end{prop}

\begin{proof} 
We must show that the section functor is fully faithful and essentially surjective. We proved essential surjectivity in Lemma \ref{lem:esssurj}. So we need to show that $\Gamma$ is surjective and injective on morphisms.

We first show that $\Gamma$ is surjective on morphisms. In other words, that every $C^\infty(M)$-linear map  $F\colon \Gamma(E)\to \Gamma(E')$ for VB-tepui $E$ and $E'$ comes from a VB-tepui morphism.  We know that $E=V/D$ and $E'=V'/D'$ for some smooth vector bundles $V$ and $V'$, and smooth singular subbundles $D$ and $D'$.  We may, without loss of generality, assume $V=M\times \mathbb R^N$, with a global frame $v_1,\dots ,v_N\in\Gamma(V)$. Since the map $\Gamma(V')\to \Gamma(E')$ is surjective, there exist sections $\sigma_1,\dots,\sigma_N\in\Gamma(V')$ with $F([v_i])=[\sigma_i]$. Extending by linearity, we obtain a morphism of the projective modules $\Gamma(V)\to \Gamma(V')$. By the classical Serre-Swan theorem, it is given by $f_*$ for a vector bundle map $f\colon V\to V'$. Hence, we obtain a diagram of modules: 
  \begin{equation*}
    \begin{tikzcd}
      {\Gamma(D)} && {\Gamma(D')} \\
      {\Gamma(V)=C^\infty(M)^N} && {\Gamma(V')} \\
      {\Gamma(E)} && {\Gamma(E')}
      \arrow["{f_*}", from=1-1, to=1-3]
      \arrow[hook, from=1-1, to=2-1]
      \arrow[hook, from=1-3, to=2-3]
      \arrow["{f_*}", from=2-1, to=2-3]
      \arrow[two heads, "{q_*}"', from=2-1, to=3-1]
      \arrow[two heads, "{q'_*}"', from=2-3, to=3-3]
      \arrow["F", from=3-1, to=3-3]
    \end{tikzcd}
  \end{equation*}
Since any element of $D$ can be extended to a section of $D$, the map $f$ must send $D$ to $D'$. Therefore, this diagram can be translated to a diagram of VB-tepui
  \begin{equation*}
    \begin{tikzcd}
      V && {V'} \\
      E && {E'}.
      \arrow["f", from=1-1, to=1-3]
      \arrow["q"', two heads, from=1-1, to=2-1]
      \arrow["{q'}"', two heads, from=1-3, to=2-3]
      \arrow["{\tilde f}", from=2-1, to=2-3]
    \end{tikzcd}
  \end{equation*}
The only thing we need to check is that $\tilde f$ is smooth, but that follows from the fact that $q'\circ f$ is smooth and $q$ is a subduction.

Now we show that $\Gamma$ is injective on morphisms. The categories of VB-tepui and $C^\infty(M)$-modules are enriched over vector spaces, and $\Gamma$ preserves the vector space structures on the Hom-sets, so it is enough to show that a non-zero morphism $f \colon E\to E'$ of VB-tepui induces a non-trivial morphism of modules of sections. The map $f$ being non-trivial means that there is $e\in E_m\subset E$ with $f(e)\neq 0_m$. On a VB-tepui, there is a section through any element, hence there is some $\sigma\in \Gamma(E)$ such that $\sigma(m)=e$. Then $f_* \sigma$ is a non-zero section of $E'$. This means that $\Gamma$ is indeed injective.

\end{proof}

\subsection{A Serre-Swan theorem for locally finitely generated modules}

We will now formulate the Serre-Swan theorem for general VB-tepui. To do so, we need to describe the class of $C^\infty(M)$-modules which can be attained as sections of a VB-tepui. Firstly, since the sections of a VB-tepui always form a sheaf, we will only consider modules which can be seen as global sections of some sheaf of modules over the sheaf of rings $C^\infty_M$ of smooth functions. Second, since VB-tepui are locally finite-dimensional, the associated sheaf of sections must be locally finitely generated (i.e., of finite-type; any point must have an open neighbourhood on which the sheaf is finitely generated). In order to carry out the translation between the category $C^\infty(M)\textrm{-}\cat{Mod}$ of  $C^\infty(M)$-modules and the category $C^\infty_M\textrm{-}\cat{Mod}$ of sheaves of $C^\infty_M$-modules, we will use the following result.

\begin{thm}\label{thm:equiv1} We consider the functors:
	\begin{itemize}
		\item $G\colon C^\infty_M\textrm{-}\cat{Mod}\to  C^\infty(M)\textrm{-}\cat{Mod}$, the global sections functor, $G(\mathcal S) \coloneqq \mathcal S(M)$
		\item  $L\colon C^\infty(M)\textrm{-}\cat{Mod} \to C^\infty_M\textrm{-}\cat{Mod}$, the localization functor, given at $Q$ by sheafifying the presheaf $U\mapsto (S_U)^{-1}Q$, where $S_U$ is the space of smooth functions on $M$ which have no zeroes on $U$.
	\end{itemize}
These functors induce an equivalence of categories between $C^\infty_M\textrm{-}\cat{Mod}$ and the essential image of $G$ in $C^\infty(M)\textrm{-}\cat{Mod}$.
\end{thm}

\begin{proof}
	This is \cite[Theorem 3.11]{gonzalessalas}, together with the observation that $C^\infty(M)$ is a differential algebra for any smooth manifold $M$ (see \cite[p.\ 37]{gonzalessalas}).
\end{proof}

We will use the following terminology.
\begin{defn} A $C^\infty(M)$-module $Q$ is called
	\begin{itemize}
		\item \define{locally finitely generated} if its localization $L(Q)$ is a locally finitely generated sheaf of $C^\infty_M$-modules;
		\item \define{global} if it is in the essential image of $G$ (equivalently, if $G(L(Q))\cong Q$).
	\end{itemize}
We will also call $G(L(Q))$ the \define{global hull} of $Q$.
\end{defn}

Given a sheaf of modules, it suffices to verify fiber-determinedness on global sections.

\begin{lem}
  \label{lem:global-geometric}
  Let $S$ be a sheaf of $C^\infty_M$-modules. If $S(M)$ is a fiber-determined $C^\infty(M)$-module, then $S(U)$ is a fiber-determined $C^\infty(U)$-module for all $U \subseteq M$.
\end{lem}

\begin{proof}
  Let $\sigma\in S(U)$, such that $\sigma\in I_m(U)S(U)$ for all $m\in U$, where $I_m(U)$ denotes the functions on $U$ that vanish at $m$. We need to show that $\sigma$ is zero, so let us assume that $\sigma\neq 0$, and derive a contradiction. First of all, $\sigma\neq 0$ means that there is a point $m_0\in U$ such that the stalk of $\sigma$ at $m_0$ is nonzero. Let us pick neighbourhoods $V\subset W\subset U$ such that $\overline V\subset W$ and $\overline W\subset U$. Let $\chi\in C^\infty(U)$ be a function with support in $\overline W$ and which is constantly $1$ on $\overline V$. The element $\sigma'\coloneqq\chi\cdot \sigma$ satisfies:
	\begin{itemize}
		\item $\sigma'\neq 0$ (since it has the same stalk as $\sigma$ in $m$)
		\item $\sigma'\in I_m(U)S(U)$ for all $m\in U$.
		\item $\sigma'|_{U\smallsetminus \overline W}=0$
	\end{itemize} 
We can glue $\sigma'$ with the function 0 on $M\smallsetminus \overline W$, to obtain a section $\sigma''\in S(M)$. A calculation shows that $\sigma''\in I_mS(M)$ for all $m$, so $\sigma''=0$ by fiber-determinedness of $S(M)$. But the stalk of $\sigma''$ at $m_0$ is the same as the stalk of $\sigma$ at $m_0$, which leads to a contradiction. Thus, $\sigma$ was zero to begin with.
\end{proof}

We can now prove our Serre-Swan theorem for general VB-tepui.
\begin{thm}\label{thm:gensingserreswan} The section functor $\Gamma$ induces an equivalence of categories between the category of VB-tepui and locally finitely generated, global, fiber-determined $C^\infty(M)$-modules.
\end{thm}
 
\begin{proof}
  The section functor $\Gamma$ factors through sheaves via the ``local section'' functor
  \begin{equation*}
  \tilde \Gamma \colon \cattwo{VB}\textrm{-}\mathrm{tepui}\to  C^\infty_M\textrm{-}\cat{Mod}.
  \end{equation*}
The essential image of $\tilde \Gamma$ is contained in the subcategory of locally finitely generated and fiber-determined sheaves of  $C^\infty_M$-modules. By Lemma \ref{lem:global-geometric} and Theorem \ref{thm:equiv1}, this is equivalent to the category of locally finitely generated, global, and fiber-determined $C^\infty(M)$-modules. So we just need to show that $\tilde \Gamma$ is an equivalence of the categories of VB-tepui and of locally finitely generated and fiber-determined sheaves of $C^\infty_M$-modules. Since morphisms of VB-tepui are fully determined by their restriction to small open sets $U\subset M$, Proposition \ref{prop:bounded-serre-swan} already implies fullness and faithfulness of the functor. Thus, we need only show essential surjectivity. Let $S$ be a locally finitely generated and fiber-determined sheaf of modules. We cover $M$ by $U_i$ such that each $S|_{U_i}$ is finitely generated. Then by Proposition \ref{prop:bounded-serre-swan}, there are VB-tepui $E_i\to U_i$, such that $S(U_i)\cong \Gamma(E_i)$. On double intersections $U_{ij}=U_i\cap U_j$, we have 
\begin{equation*}
  \Gamma((E_i)_{U_{ij}})\cong S(U_{ij})\cong \Gamma((E_j)_{U_{ij}}),
\end{equation*}
which induce isomorphisms $\phi_{ij}\colon (E_i)_{U_{ij}}\to (E_j)_{U_{ij}}$. We now define $E\coloneqq \left(\coprod_i E_i\right)/{\sim}$, where the equivalence is generated by $\phi_{ij}$. The resulting object is a VB-tepui with the desired property.

\end{proof}

\subsection{The smooth section functor and the topological tensor product}

In this subsection, we address two possible extensions of our Serre-Swan theorem. First, we observe that, given a diffeological bundle of vector spaces $E \to M$, we have treated $\Gamma(E)$ purely algebraically. However, since $\cat{Dflg}$ is Cartesian closed, the space $\Gamma(E)$ can also carry its natural functional diffeology. As we saw in Proposition \ref{prop:frechet-vs-diffeology}, this coincides with its usual Fr\'{e}chet structure when $E = V$ is a vector bundle, and the result is that $\Gamma(V)$ is a Fr\'{e}chet module over the Fr\'{e}chet algebra $C^\infty(M)$. We will show that, more generally, $\Gamma(E)$ is a Fr\'{e}chet module when $E$ is a VB-tepui.

Second, for the classical Serre-Swan theorem, one also has that the relevant functor $\Gamma$ is monoidal. For VB-tepui, this is not the case for the algebraic tensor product of modules, nor for the completed tensor product of Fr\'{e}chet modules; we give an example illustrating this.

\begin{lem}
  \label{lem:lifting-sections}
  Let $V \to M$ be a vector bundle, and $D \subseteq V$ be a smooth singular subbundle. There is a natural diffeomorphism and isomorphism of $C^\infty(M)$-modules $I\colon \Gamma(V/D) \to \Gamma(V)/\Gamma(D)$ fitting in the diagram below.
   \begin{equation*}
    \begin{tikzcd}
      & \Gamma(V) \ar[dl, "q_*"'] \ar[dr, two heads] & \\
      \Gamma(V/D) \ar[rr, "I"] & & \Gamma(V)/\Gamma(D).
    \end{tikzcd}
  \end{equation*}
Here $q\colon V \to V/D$ is the quotient map.
\end{lem}

\begin{proof}
The map $I$ is a well-defined bijection because two sections of $V$ agree when projected to $V/D$ if and only if their difference is a section of $D$. If $q_*$ is a subduction, then the universal property for subductions implies that $I$ is a diffeomorphism. It is an isomorphism of $C^\infty(M)$-modules because (assuming that $q_*$ is a subduction) the downward arrows are epimorphisms in $C^\infty(M)\mathrm{-}\cat{Mod}$. 

We show that $q_*\colon \Gamma(V) \to \Gamma(V/D)$ is a subduction. Fix a plot $\sigma\colon U \to \Gamma(V/D)$, and denote its transpose as usual by $\hat{\sigma}$.  This is smooth, and since $q\colon V \to V/D$ is a subduction, there is a locally finite cover $\{U_i \times M_i\}$ of $U \times M$, equipped with lifts $\lambda_i\colon U_i \times M_i \to V$ of $\hat{\sigma}$ along $q$. Take a partition of unity $\{\chi_i\}$ subordinate to $\{U_i \times M_i\}$. Then the map
  \begin{equation*}
    U \times M \to V, \quad (r,m) \mapsto \sum_i \chi_i(r,m)\lambda_i(r,m)
  \end{equation*}
  is smooth, and the associated map $U \to \Gamma(V)$ is the desired lift of $\sigma$.
\end{proof}

\begin{prop}\label{prop:frechetmodfromsections}
  Let $E \to M$ be a VB-tepui. Then $\Gamma(E)$, equipped with the D-topology underlying its standard functional diffeology, is a Fr\'{e}chet module. Therefore, taking sections gives the functor
  \begin{equation*}
    \Gamma \colon \cattwo{VB}\textrm{-tepui} \to C^\infty(M)\textrm{-}\cat{Frch}\textrm{-}\cat{Mod}.
  \end{equation*}
  Moreover, the essential image of $\Gamma$ consists of \emph{nuclear} Fr\'{e}chet spaces.
\end{prop}

We include the adjective ``nuclear'' for those who are familiar with the term, and will not expand the definition here. The reader may refer to \cite[Chapter 50]{Trev67} for details on nuclear spaces.

\begin{proof}
We first assume that $E \to M$ has globally bounded rank. By Corollary \ref{cor:globfin}, there is a smooth vector bundle $V \to M$ and a smooth singular subbundle $D$ of $V$ such that $E \cong V/D$. By Lemma \ref{lem:lifting-sections} above, we have $\Gamma(V/D) \cong \Gamma(V)/\Gamma(D)$ as diffeological spaces. Therefore, the same is true for their underlying topological spaces. The underlying space for $\Gamma(V)$ is its usual Fr\'{e}chet structure by Lemma \ref{prop:frechet-vs-diffeology} and Remark \ref{rk:d-topology-is-mackey}, and this structure is nuclear (\cite[Corollary to Theorem 51.4]{Trev67}). The space $\Gamma(D)$ is a closed subset of $\Gamma(V)$, since it is the intersection of the sets $\{\sigma\in\Gamma(V)\mid \sigma(m)\in D_m\}$, each of which is closed as the inverse image of the closed subset $D_m\subset E_m$ via the (continuous) evaluation map $\operatorname{ev}_m\colon \Gamma(E)\to E_m$. Therefore, $\Gamma(V)/\Gamma(D)$ is again a nuclear Fr\'{e}chet space (\cite[Proposition 50.1]{Trev67}), hence so is $\Gamma(V/D)$. 

Now, locally any VB-tepui $E$ is a quotient $V/D$. Choose a countable cover $\{U_i\}$ of $M$ by open subsets, such that each $E|_{U_i}$ arises as a quotient bundle. The projective limit of the $\Gamma(E|_{U_i})$ in $\cat{Dflg}$ coincides with the standard functional diffeology, because to be a smooth section is a local condition. On the other hand, by the recollement theorem \cite[Appendix A.1]{gonzalessalas}, we can also equip $\Gamma(E)$ with the nuclear Fr\'{e}chet structure induced by taking the projective limit of the $\Gamma(E|_{U_i})$ in $\cat{Frch}$ (cf.\ \cite[Proposition 50.1]{Trev67}). Because the functor $\sigma_\infty\colon \cat{Frch} \to \cat{Dflg}$ is a right adjoint (by Lemma \ref{lem:adjunctions}), the diffeology associated to this Fr\'{e}chet projective limit structure coincides with the functional diffeology. By Remark \ref{rk:d-topology-is-mackey}, we see that the D-topology of $\Gamma(E)$ coincides with the Fr\'{e}chet topology of $\Gamma(E)$.  
\end{proof}

\begin{rk}
  \label{rk:coefficient-diffeology}
  In Proposition \ref{prop:frechetmodfromsections}, we saw two equivalent descriptions of the diffeology on $\Gamma(E)$ for a VB-tepui $E \to M$. Here we give yet one more. Let $p\colon W \to \Gamma(E)$ be a plot, and fix $m_0 \in M$. Let $U$ be a neighbourhood of $m_0$ in $M$ over which the restricted bundle $E_U \to U$ has the form $E_U \cong V/D$ for a vector bundle $V \to U$ and a smooth singular subbundle $D \subseteq V$. Then, near each point of $W$, we can lift $p|_U$ (the restriction of $p$ to $U$) as in the following diagram.
  \begin{equation}
    \label{eq:liftability}
    \begin{tikzcd}
      & \Gamma(V) \ar[dr, "q_*"] & \\
      & & \Gamma(V/D) \ar[d, "\cong"] \\
      W \ar[uur, dashed]\ar[r, "p"] & \Gamma(E) \ar[r, "|_U"] & \Gamma(E_U).
    \end{tikzcd}
  \end{equation}
  The lift is possible because $p|_U$ is a plot of $\Gamma(E_U)$, and by Lemma \ref{lem:lifting-sections} the map $q_*$ is a subduction.

  The collection of plots of $\Gamma(E)$ for which the lifting problem \eqref{eq:liftability} has a solution is called the \emph{coefficient diffeology} in \cite{GV22}. There, the authors define the coefficient diffeology for a general $C^\infty(M)$-module, and we have justified that in the case of a VB-tepui, the natural diffeology on $\Gamma(E)$ coincides with its coefficient diffeology.
\end{rk}
\begin{rk} If we describe the singular subbundle $D$ as the image $r(W)$ of a vector bundle morphism $W\to V$, then $\Gamma(V/D)=\Gamma(V/r(W))$ is in general not equal to $\Gamma(V)/r_*\Gamma(W)$. The latter might not be fiber-determined, e.g., if one takes $M=\mathbb R$, with $V=W=M\times \mathbb R$ and $r(x,v)=(x,x^2v)$.
\end{rk}

We end this section by discussing how the functor $\Gamma$ interacts with tensor products. Let us start by recalling the tensor product of diffeological bundles of vector spaces, as described in \cite[Section 3.2.6]{ChrWu23}.

\begin{defn}
  \label{def:tepui-tensor}
  Given two diffeological bundles of vector spaces $E$ and $E'$ over a manifold $M$, their \define{tensor product} $E \otimes E'$ has underlying set $|E \otimes E'| \coloneqq \coprod |E_m \otimes_{\R} E_m'|$, where the tensor product of the spaces in the disjoint union is the algebraic tensor product of vector spaces. The plots are the maps that are locally of the form
  \begin{equation*}
    U \to E \otimes E', \quad u \mapsto \sum_i r_i(u) q_i(u) \otimes q_i'(u), 
  \end{equation*}
  for smooth maps $r_i\colon U \to \R$ and $q\colon U \to E$ and $q'\colon U \to E'$ such that $\pi \circ q = \pi' \circ q'$.
\end{defn}

\begin{lem}
  \label{lem:tensor-is-tepui}
  If $E$ and $E'$ are VB-tepui over $M$, then $E \otimes E'$ is also VB-tepui.
\end{lem}
\begin{proof}
  Since to be a VB-tepui is a local condition, we may use Lemma \ref{cor:globfin} and assume that $E \cong V/D$ and $E' \cong V'/D'$ for smooth vector bundles $V$ and $V'$ over $M$ and smooth singular subbundles $D\subseteq V$ and $D'\subseteq V'$. Working fiberwise, and using the fact that $|E_m \otimes E_m'|$ is the algebraic tensor product, yields the horizontal bijection in the diagram
  \begin{equation}
    \label{eq:quotient-tensor-bundle-iso}
    \begin{tikzcd}
      & V\otimes V' \ar[dl, "\mathrm{quot}\otimes\mathrm{quot}"'] \ar[dr, two heads] & \\
      E \otimes E' \ar[rr, "\cong"] & & V\otimes V' / \langle V \otimes D' + D \otimes V' \rangle.
    \end{tikzcd}
  \end{equation}
The tensor product of subductions remains a subduction, so we see that our bijection is a diffeomorphism. Since the singular distribution $V \otimes D' + D \otimes V' \subseteq V \otimes V'$ is smooth, we see that $E \otimes E'$ is VB-tepui.
\end{proof}

Since the external tensor product $E\boxtimes E'\to M\times M'$ of two bundles $E\to M$ and $E'\to M'$ can be defined as $\pi_1^*E\otimes \pi_2^*E'$ for the projections $\pi_1\colon M\times M'\to M$ and $\pi_2\colon M\times M'\to M'$, together with Lemma \ref{lem:pullbackvb}, the above implies:

\begin{cor}
If $E$ and $E'$ are VB-tepui over $M$, respectively $M'$, then $E \boxtimes E'$ is also a VB-tepui.
\end{cor}

We can now show that the functor $\Gamma$ is not monoidal, where we use the tensor product of VB-tepui on the domain and the tensor product of $C^\infty(M)$-modules on the codomain.

\begin{ex}\label{ex:not-monoidal}
  Let $V \to \R$ be the trivial $\R$-bundle over $\R$. Fix two smooth singular distributions
  \begin{equation*}
    D_\geq \coloneqq \{(x,v) \mid \text{ if } x \geq 0 \text{ then } v = 0\} \quad D_{\leq} \coloneqq \{(x,v) \mid \text{ if } x \leq 0 \text{ then } v=0\}.
  \end{equation*}
So, sections of $D_{\geq}$ vanish on $\{x \geq 0\}$, and sections of $D_{\leq}$ vanish on $\{x \leq 0\}$.  For the VB-tepui $E_{\geq} \coloneqq V/D_{\geq}$ and $E_{\leq} \coloneqq V/D_{\leq}$, we claim that
  \begin{equation*}
    \Gamma(E_\geq) \mathbin{{\otimes}} \Gamma(E_\leq) \cong \R^{\mathbb{N}} \quad \text{whereas} \quad \Gamma(E_\geq \otimes E_\leq) \cong \R.
  \end{equation*}
  This entails two computations. We denote by $A \coloneqq \Gamma(V) \cong C^\infty(\R)$, and note that all tensor products in this example are over $A$. 
  
  For the first computation, we start with the tensor product of the modules:
  \begin{align*}
    \Gamma(E_\geq) \mathbin{{\otimes}} \Gamma(E_\leq)  &\underset{(1)} \cong (A/\Gamma(D_\geq)) \mathbin{{\otimes}} (A/\Gamma(D_\leq)) \\
    &\underset{(2)}\cong \frac{A \mathbin{{\otimes}} A}{{\langle A \mathbin{{\otimes}} \Gamma(D_\leq) + \Gamma(D_\geq) \mathbin{{\otimes}} A\rangle}}\\
    &\underset{(3)}\cong \frac{A}{{\langle \Gamma(D_\leq) + \Gamma(D_\geq)\rangle}} \\
    &\underset{(4)}\cong \frac{A}{\{f \mid f^{(n)}(0) = 0\}} \\ 
    &\underset{(5)}\cong \R^{\mathbb{N}}.
  \end{align*}
  Isomorphism (1) follows from Lemma \ref{lem:lifting-sections},  (2) is a general property of the tensor product, and (3) is the neutrality of $A$ with respect to $A$. Isomorphism (4) follows from the fact that any function which is flat (i.e., vanishes with all its derivatives) at zero, can be written in a unique way as the sum of a function vanishing on $(-\infty,0]$ with a function vanishing on $[0,\infty)$. Isomorphism (5) reflects the fact that any Taylor expansion at a point is realized by some smooth function (Borel's theorem).
  
 For the VB-tepui, we have:
  \begin{align*}
    \Gamma(E_\geq \otimes E_\leq) &\cong \Gamma \left(\frac{V \otimes V}{V \otimes D_\leq + D_\geq \otimes V}\right) &&\text{by \eqref{eq:quotient-tensor-bundle-iso} in Lemma \ref{lem:tensor-is-tepui}} \\
    &\cong \Gamma\left(\frac{V}{D_\leq + D_\geq}\right) \\
  &\cong \frac{A}{\{f \mid f(0) = 0\}} &&\text{by Lemma \ref{lem:lifting-sections}} \\
    &\cong \R.
  \end{align*}
\end{ex}

\begin{rk} This example shows that the tensor product of two fiber-determined modules need not be fiber-determined, thus providing a counterexample to \cite[Exercise on tensor products, page 178]{nestruev}.
\end{rk}

The above example illustrates that the presence of invisible elements may obstruct the monoidality of $\Gamma$. In fact, this turns out to be essentially the only obstruction.

\begin{prop}\label{prop:geometrization} Let $E\to M$ and $E'\to M$ be VB-tepui and $\alpha\colon \Gamma(E)\otimes_{C^\infty(M)} \Gamma(E')\to \Gamma(E\otimes E')$ be the natural assignment $\sigma_1\otimes \sigma_2\mapsto \sigma_1\sigma_2$. Then $\alpha$ induces an isomorphism from the global hull of the fiber-determination of  $\Gamma(E)\otimes_{C^\infty(M)} \Gamma(E')$ to $\Gamma(E\otimes E')$.
\end{prop}

\begin{proof} We first treat the case of bounded rank VB-tepui. We start by picking vector bundles $V$ and $V'$ and singular distributions $D$ and $D'$ such that $E\cong V/D$ and $E'\cong V'/D'$. Then, using $\Gamma(V)\otimes \Gamma(V') \cong \Gamma(V\otimes V')$, by e.g. \cite[Theorem 12.39]{nestruev}, we obtain 
\begin{align*}
\Gamma(E)\otimes_{C^\infty(M)} \Gamma(E')&=\frac{\Gamma(V\otimes V')}{\Gamma(V)\otimes \Gamma(D')+\Gamma(D)\otimes \Gamma(V')}\\
\Gamma(E\otimes E')&=\frac{\Gamma(V\otimes V')}{\Gamma(\langle V\otimes D' + D\otimes V'\rangle)}
\end{align*}
This already implies that $\alpha$ is surjective. Moreover, since $\Gamma(E\otimes E')$ is fiber-determined, all invisible elements in $\Gamma(E)\otimes_{C^\infty(M)} \Gamma(E')$ lie in the kernel of $\alpha$. So the only thing we have to show is that all elements in the kernel of $\alpha$ are invisible, which in turn reduces to the fact that the maps 
\begin{equation*}
  i_m^*\alpha\colon i_m^*(\Gamma(E)\otimes_{C^\infty(M)} \Gamma(E'))\to i_m^*\Gamma(E\otimes E')
\end{equation*}
are injective for all $m\in M$. Consider the following diagram:
\[\begin{tikzcd}
	{\Gamma(V)\otimes \Gamma(D')\oplus \Gamma(D)\otimes \Gamma(V')} & {\Gamma(V\otimes V')} & {\Gamma(E)\otimes_{C^\infty(M)}\Gamma(E')} \\
	{i_m^*\left(\Gamma(V)\otimes \Gamma(D')\oplus\Gamma(D)\otimes \Gamma(V')\right)} & {i_m^*(\Gamma(V\otimes V'))} & {i_m^*\left(\Gamma(E)\otimes_{C^\infty(M)}\Gamma(E')\right)} \\
	{V_m\otimes D'_m+D_m\otimes V'_m} & {V_m\otimes V'_m} & {\frac{V_m\otimes V'_m}{V_m\otimes D'_m+D_m\otimes V'_m}}
	\arrow[from=1-1, to=1-2]
	\arrow[two heads, from=1-1, to=2-1]
	\arrow[two heads, from=1-2, to=1-3]
	\arrow[two heads, from=1-2, to=2-2]
	\arrow[two heads, from=1-3, to=2-3]
	\arrow[from=2-1, to=2-2]
	\arrow[two heads, from=2-1, to=3-1]
	\arrow[two heads, from=2-2, to=2-3]
	\arrow[equals, from=2-2, to=3-2]
	\arrow["{i_m^*\alpha}", from=2-3, to=3-3]
	\arrow[from=3-1, to=3-2]
	\arrow[two heads, from=3-2, to=3-3]
\end{tikzcd}\]
Here, the horizontal rows are quotients, the second row is obtained from the first by applying $i_m^*$, and the bottom row is obtained by evaluating at $m$. For the surjectivity of the bottom arrow in the left column, we use the fact that $D$ and $D'$ are \emph{smooth} singular subbundles. Now, a diagram chase shows  that $i_m^*\alpha$ is bijective, implying that the kernel of $\alpha$ consists exactly of the invisible elements. 

In the case of unbounded rank, the above still holds locally.  On $U$ a sufficiently small open set, $\Gamma((E\otimes E')|_U)$ is isomorphic to the fiber-determination of $\Gamma(E|_U)\otimes_{C^\infty(U)}\Gamma(E'_U)$. We obtain the global statement by applying the global hull functor to both sides of the isomorphy and observing that $\Gamma(E\otimes E')$ is already global.
\end{proof}

As we have seen in Proposition \ref{prop:frechetmodfromsections}, for a VB-tepui $E$, the space of sections $\Gamma(E)$ is not just a $C^\infty(M)$-module but a Fréchet-module. In particular, we could compare the tensor product of VB-tepui, not only with the usual tensor product $\Gamma(E)\otimes_{C^\infty(M)}\Gamma(E')$ of $C^\infty(M)$-modules, but also to the completed tensor product $\Gamma(E) \hat\otimes_{C^\infty(M)} \Gamma(E')$ of Fréchet modules (cf.\ e.g.\ \cite[Section 6.2]{gonzalessalas}).

\begin{defn}
	\label{def:completed-module-tensor}
	Given two $C^\infty(M)$ locally convex modules $Q$ and $Q'$, their \define{completed projective tensor product} $Q \mathbin{\hat{\otimes}}_{C^\infty(M)} Q'$ is the completion of the locally convex topological vector space $Q \otimes_{C^\infty(M)} Q'$, which is the algebraic tensor product equipped with the final topology induced by the canonical map $Q \otimes_\R Q'\to Q \otimes_{C^\infty(M)} Q'$, whose domain carries the finest locally convex topology such that the canonical map $Q \times Q' \to Q\otimes_\R Q'$ is continuous.
\end{defn}
The completed projective tensor product of two Fr\'{e}chet $C^\infty(M)$-modules is again Fr\'{e}chet, so we have a tensor product in the category of Fr\'{e}chet $C^\infty(M)$-modules.

\begin{prop} In the setting of Proposition \ref{prop:geometrization}, the map $\alpha\colon \Gamma(E)\otimes_{C^\infty(M)} \Gamma(E')\to \Gamma(E\otimes E')$ factors through the (surjective) completion map $c\colon \Gamma(E)\otimes_{C^\infty(M)} \Gamma(E')\to \Gamma(E)\hat\otimes_{C^\infty(M)} \Gamma(E')$, i.e., there is a homomorphism $\beta\colon \Gamma(E)\hat\otimes_{C^\infty(M)} \Gamma(E')\to \Gamma(E\otimes E')$, such that $\alpha=\beta\circ c$. In particular, $\Gamma(E\otimes E')$ is the global hull of the fiber-determination of $\Gamma(E)\hat\otimes_{C^\infty(M)} \Gamma(E')$.
\end{prop}

\begin{proof} We will again prove the assertion in the bounded rank case, the general statement following by a sheaf theoretic argument. So let $E \cong V/D$ and $E' \cong V'/D'$. We first note that $\Gamma(V)\hat\otimes_{C^\infty(M)} \Gamma(V')=\Gamma(V\otimes V')$, cf. e.g. \cite[Theorem 2.3.8.]{herscovichRenormalizationQuantumField2019}. In particular, this means that 
\begin{align*}
	{\Gamma(E)\hat\otimes_{C^\infty(M)} \Gamma(E')}={\left(\frac{\Gamma(V\otimes V')}{Q}\right)}^\wedge
= \frac{{\Gamma(V\otimes V')}^\wedge}{{Q}^\wedge}=\frac{{\Gamma(V\otimes V')}}{\overline{Q}}
\end{align*}
where we write used the notation $Q=\Gamma(V)\otimes \Gamma(D')+\Gamma(D)\otimes \Gamma(V')$, seen as a subspace of ${\Gamma(V\otimes V')}$, and denoted completions by the ``$~^\wedge~$'' superscript. In particular, this means that the map $c$ is a quotient map. For the map $\alpha$ to factor through the completion, we only need to show that the kernel of the completion map $c$ is contained in the set of invisible elements of $\Gamma(E)\otimes_{C^\infty(M)} \Gamma(E')$. By \cite[Theorem 1.3, Chapter V]{tougeronIdeauxFonctionsDifferentiables1972}\footnote{There the theorem is stated for submodules of trivial modules, but since Q is a submodule of a projective module (which is a direct summand in a free one)  the assertion of the theorem still holds in our case}, the closure $\bar Q$ of $Q$  consists of elements $\sigma$ in $\Gamma(V\otimes V')$, whose infinite jet $j^\infty_m\sigma$ at $m$ coincides with the $j^\infty_mq$ for some $q\in Q$. In particular, $\sigma(m)=q(m)$. But $q(m)\in \ker(V_m\otimes V_{m}'\to E_m\otimes E_m')$, i.e., $[\sigma(m)]=0\in E_m\otimes E_m'$ for all $m$. This means that $\sigma$ is invisible and concludes the proof.
\end{proof}

\begin{rk} In Example \ref{ex:not-monoidal}, the denominator $Q=\{f \mid f^{(n)}=0\}$ was already closed, in particular it provides an example where $\Gamma(E)\hat \otimes_{C^\infty(M)}\Gamma(E')$ is equal to $\Gamma(E)\otimes_{C^\infty(M)}\Gamma(E')$. Thus, both have invisible elements and are hence not isomorphic to $\Gamma(E\otimes E')$. 
\end{rk}

\subsection{Algebraic and geometric base-change}

Given a bundle of vector spaces $E \to X$ over a diffeological space $X$ and a smooth map $\pi$ from some other space $P$ to $X$, there are two ways to obtain a module over $C^\infty(P)$: either pulling back $E \to X$ along $\pi$ and then applying the global sections functor, or applying the algebraic base change $\pi^*\Gamma(E) \coloneqq \gi(P)\otimes_{\gi(X)} \Gamma(E)$ directly to the section space $\Gamma(E)$. These two are related by the natural map 
\begin{align}\label{eq:defalpha}
\alpha_E\colon \pi^*\Gamma(E) \to \Gamma(\pi^* E), \quad \sum_if_i\otimes \sigma_i\mapsto \sum_if_i\cdot \pi^*\sigma_i.
\end{align}

For vector bundles over smooth manifolds, this map is an isomorphism \cite[Theorem 12.43]{nestruev}. However, for general bundles of vector spaces, and even VB-tepui, $\alpha$ is neither injective nor surjective.
\begin{ex}[Non-surjectivity of $\alpha$]
Let $\R^\infty$ denote the diffeological vector space of real sequences with finitely many non-zero terms. Let $\{e_n\}_{n \in \mathbb{N}}$ denote the standard basis of $\R^\infty$. Fix smooth functions $\{f_n\colon \R \to \R\}$ so that each has $n$-zeros, namely $f_n^{-1}(0) = \{i \in \mathbb{N} \mid i \leq n\}$. We form the singular subbundle $D$ of the trivial bundle $\R \times \R^\infty$ by declaring $D_x \coloneqq \{x\} \times \operatorname{span} \{f_n(x)e_n \mid n \in \mathbb{N}\}$. Consider the VB-tepui $E \coloneqq (\R \times \R^\infty)/D$ over $M \coloneqq \R$. We have
\begin{equation*}
  E_x\cong \begin{cases}
    \mathbb R^n & \text{if } x=n\in\mathbb N, \\
    0  & \text{else}.
  \end{cases}
\end{equation*}
Let $P\coloneqq \mathbb R$ and $\pi\colon P\to M$ be a smooth surjective map which passes each integer point $n\in\mathbb N$ exactly $n$ distinct times, say $t^n_1< \dots <t^n_n$. For concreteness, we may take $\pi$ as follows: let $T_n = \frac{n(n+1)}{2}$ denote the $n$-th triangular number, and say that $\pi(x) = n \in \mathbb{N}$ if and only if $x$ is natural and $T_{n-1} < x \leq T_n$. For each $i$, let $h_i\colon \R \to [0,1]$ be a smooth function supported on $\bigcup_{n=1}^\infty (t^n_i - 1/2, t_i^n+1/2)$, where $h_i(x) = 1$ if and only if $x = t_i^n$. Set $\bar{\sigma}$ to be the section of $\pi^*E$ given by $\bar{\sigma}(z) \coloneqq (z,\pi(z), \sum_i h_i(z)[e_i])$.

We claim that the element $\bar{\sigma}$ can never be attained as $\alpha(\sum_{i=1}^kf_i\otimes \sigma_i)$. To see this, we identify sections of $\pi^*E$ with functions $\sigma \colon \R \to E$, and then consider the set

\begin{equation*}
 \bigcup_{x \in \R} \operatorname{span}\{\sigma(z) \mid \pi(z)=x\} \subset E
\end{equation*}

For $\sigma=\bar\sigma$, this set is equal to $E$ at any $x$, yet for any $\sigma=\alpha(\sum_{i=1}^kf_i\otimes \sigma_i)$, is it at most $k$-dimensional at any point $x$, hence not $E$ since $E$ is not of bounded rank.
\end{ex}

\begin{ex}[Non-injectivity of $\alpha$]\label{ex:algpullbnotgeom}
  We use an example similar to Example \ref{ex:not-monoidal}. Take $V \to M$ to be the trivial $\R$-bundle over $\R$, and fix the smooth singular distribution
  \begin{equation*}
    D \coloneqq D_{\leq} = \{(x,v) \mid \text{ if } x \leq 0 \text{ then } v=0\}.
  \end{equation*}
  Let $P \coloneqq \R$, and take $\pi\colon P \to M$ to be the map $x \mapsto x^2$. We now pull back the VB-tepui $E \coloneqq V/D$ along  $\pi$. We first observe that the fiber of $\pi^*D$ is $\R$ over $x \neq 0$, and is trivial over $x=0$, thus $\pi^*E \cong \frac{\pi^*V}{\pi^*D}$ has only one non-trivial fiber: the one over $x = 0$.

The natural map $C^\infty(P)\otimes_{C^{\infty}(M)}\Gamma(D)\to \Gamma(\pi^*D)$ is not surjective. For an element $\sum f_i\otimes \sigma_i$ in the domain, all the $\sigma_i$ have vanishing $\infty$-jet at the origin, hence the image $\sum f_i \cdot \pi^*\sigma_i$ will also have vanishing $\infty$-jet at the origin. Conversely, the section $f$ of $\Gamma(\pi^*D)$ corresponding to the identity $x \mapsto x$ vanishes linearly at the origin.
\end{ex}

For VB-tepui, the non-surjectivity of $\alpha$ is a minor technical issue that will be resolved by passing to the global hull of $ \gi(P)\otimes_{\gi(M)} \Gamma(E)$, so let us deal with the non-injectivity of $\alpha$. It is, in general, due to the fact that $ \gi(P)\otimes_{\gi(M)} \Gamma(E)$ may not be fiber-determined. This is resolved by passing to its fiber-determination.

\begin{prop}\label{prop:basechange} Let $E \to M$ be a VB-tepui of globally bounded rank and $\pi\colon P \to M$ be a smooth map. Then the map $\alpha$ defined in \eqref{eq:defalpha} is surjective, and $\ker(\alpha) = \mathrm{inv}(\pi^*\Gamma(E))$. In other words, $\Gamma(\pi^*E)$ is the fiber-determination of $\pi^*\Gamma(E)$.
\end{prop}

\begin{proof}
	By Corollary \ref{cor:globfin}, there exists a vector bundle $V \fto M$ and a smooth singular subbundle $D\subset V$ such that $E\cong V/D$. Then:
	\begin{equation*}
		C^\infty(P)\otimes_{C^{\infty}(M)}\Gamma(E) \cong C^\infty(P)\otimes_{C^{\infty}(M)}\frac{\Gamma(V)}{\Gamma(D)} \cong \frac{C^\infty(P)\otimes_{C^{\infty}(M)}\Gamma(V)}{C^\infty(P)\otimes_{C^{\infty}(M)}\Gamma(D)} \cong \frac{\Gamma(\pi^*V)}{\alpha_V(\, C^\infty(P)\otimes_{C^{\infty}(M)}\Gamma(D)\, )}. 
	\end{equation*}
        Here we view $\Gamma(D)$ as a submodule of $\Gamma(V)$, and $\alpha_V$ denotes the natural map from Equation \eqref{eq:defalpha} with respect to the vector bundle $V$. Inside $\Gamma(\pi^*V)$, it is clear that $\alpha_V(\, C^\infty(P)\otimes_{C^{\infty}(M)}\Gamma(D)\,)\subset \Gamma(\pi^*D)$, and therefore the map
  \begin{equation*}
    \begin{tikzcd}
      C^\infty(P)\otimes_{C^{\infty}(M)}\Gamma(E)\ar[r, "\cong"] & \frac{\Gamma(\pi^*(V))}{\alpha_V(C^\infty(P) \otimes_{C^\infty(M)} \Gamma(D))} \ar[r, two heads] &  \frac{\gs(\pi^* V)}{\gs(\pi^*D)} \cong \gs(\pi^*E),
    \end{tikzcd}
  \end{equation*}
which coincides with $\alpha_E$, is surjective.

For the second assertion, we observe that any invisible element $\sigma\in \mathrm{inv}(\pi^*\Gamma(E))$ lies in the kernel of $\alpha$. This is because any module homomorphism maps invisible elements to invisible elements, and $\Gamma(\pi^*E)$ is fiber-determined. What remains to show is that any element in the kernel of $\alpha$ is invisible. To see that, we observe that for $C^\infty(P)$-modules $Q$, the quotient $\frac{Q}{I_pQ}$ for $p\in P$ can also be written as $i_p^*Q$, where $i_p\colon \{p\}\to P$ is the inclusion.
\begin{align*}
	&\text{On the one hand we have } i_p^*\pi^*\Gamma(E)\cong i_{\pi(p)}^*\Gamma(E)=E_{\pi(p)}, \text{ because } E \text{ is VB-tepui.} \\
	&\text{On the other we have }  i_p^*\Gamma(\pi^*E)=(\pi^*E)_p=E_{\pi(p)}, \text{ because } \pi^*E \text{ is VB-tepui.}
\end{align*}
Therefore, $i_p^*\alpha$ is the identity map. Now let $\sigma\in\ker(\alpha)$, and thus $i_p^*\sigma\in \ker(i_p^*\alpha)$ for all $p$. But since $i_p^*\alpha$ is injective, that means $i_p^*\sigma=0$ for all $p$, hence $\sigma$ is invisible. 

\end{proof}

      As in Proposition \ref{prop:geometrization}, the general case follows from the bounded rank one, upon applying the global hull functor:
\begin{cor}\label{cor:basechange} Let $E\to M$ be any VB-tepui and $\pi\colon P\to M$ any smooth map. Then $\Gamma(\pi^*E)$ is the global hull of the fiber-determination of $\pi^*\Gamma(E)$.
\end{cor}

The next lemma, which derives from the proof of Proposition \ref{prop:geometrization}, gives us a tool to check whether $\alpha_E$ is an isomorphism.

\begin{lem}
  \label{lem:coker}
  Let $V \to M$ be a vector bundle, $D \subset V$ be a smooth singular subbundle, and $\pi\colon P \to M$ be a smooth map. Then the kernel of $\alpha_{V/D}$ is isomorphic to the cokernel of $\alpha_D\colon \pi^*\Gamma(D) \to \Gamma(\pi^*D)$. Therefore, $\alpha_{V/D}$ is an isomorphism if and only if $\alpha_D$ is surjective.
\end{lem}

\begin{proof}
  This follows immediately from a general fact about modules: if $M$ and $N$ are modules over a ring $R$, and $A$ and $B$ are submodules of $M$ and $N$, respectively, and $f\colon M \to N$ is a morphism with $f(A) \subseteq B$, then denoting $\bar{f}\colon M/A \to N/B$ the induced map, we have a morphism
  \begin{equation*}
    \ker(\bar{f}) \to B/f(A) = \textrm{coker}(f\colon A \to B), \quad m+A \mapsto f(m) + f(A).
  \end{equation*}
  This is an isomorphism if $f$ is an isomorphism. In the setting of this lemma, $f = \alpha_V$, which is an isomorphism, and $\bar{f} = \alpha_{V/D}$, and the restriction $\alpha_V\colon \pi^*\Gamma(D) \to \Gamma(\pi^*D)$ coincides with $\alpha_D$.
\end{proof}
We now turn to a concrete case where $\pi^*\Gamma(E)$ and $\Gamma(\pi^*E)$ coincide.

\begin{prop}\label{prop:welldefinepullback} Let $V\fto M$ be a vector bundle and $D\subset V$ a smooth singular subbundle. If $\pi \colon P \to M$ is a submersion, and if $\Gamma(D)$ is finitely generated, then $\alpha_{V/D}$ is an isomorphism, i.e.,
	\begin{equation*}
		\pi^*\Gamma(V/D) = C^\infty(P)\otimes_{C^{\infty}(M)}\Gamma(V/D)\cong \gs(\pi^*V/\pi^*D). 
	\end{equation*}	
\end{prop}

\begin{proof} By Proposition \ref{prop:basechange} we know that $\alpha_{V/D}$ is surjective, so we only need to show it is injective. By Lemma \ref{lem:coker}, it suffices to show that $\alpha_D$ is surjective. As a first step, we assume that $P = M \times L$ for some manifold $L$, and that $\pi\colon M \times L \to M$ is the projection. Since $D\subset V$, $\gs(D)$ is a global and fiber determined $\gi(M)$ module. By assumption, $\Gamma(D)$ is also finitely generated, hence by Proposition \ref{prop:bounded-serre-swan} we have a VB-tepui $E_0 \to M$ of globally bounded rank and an isomorphism of modules $\varphi \colon \Gamma(E_0) \to \Gamma(D)$. We consider the following diagram
    \begin{equation}\label{eq:diagramsurge}
    \begin{tikzcd}
      \pi^*\Gamma(E_0) \ar[r, "\pi^*\varphi"] \ar[d, "\alpha_{E_0}"] & \pi^*\Gamma(D) \ar[d, "\alpha_D"] \\
      \Gamma(\pi^*E_0) \ar[d, "\beta_{E_0}"] & \Gamma(\pi^*D) \ar[d, "\beta_D"] \\
      C^\infty(L, \Gamma(E_0)) \ar[r, "\varphi_*"] & C^\infty(L, \Gamma(D)),
    \end{tikzcd}
  \end{equation}
where $\beta_F$ denotes, for any diffeological bundle of vector spaces $F \to M$, the natural map
  \begin{equation*}
    \begin{tikzcd}
      \beta_F\colon \Gamma(\pi^*F) \ar[r] & C^\infty(L, \Gamma(F)).
    \end{tikzcd}
  \end{equation*}
By Lemma \ref{lem:coker}, it suffices to show that $\alpha_D$ is surjective. Proposition \ref{prop:basechange} implies that $\alpha_{E_0}$ is surjective, hence we can show the surjectivity of $\alpha_D$ by showing that every other map in Diagram \eqref{eq:diagramsurge} is an isomorphism. One can verify that $\beta_F$ is a diffeomorphism for any bundle of vector spaces $F\to M$, when its domain and codomain carry their functional diffeologies. Moreover, $\pi^*\varphi$ is an isomorphism because $\varphi$ is an isomorphism. It remains to justify that $\varphi_*$ is a well-defined isomorphism. The map $\varphi$ is smooth because, denoting $\iota\colon D \to V$ the inclusion, the composite map
  \begin{equation*}
    \begin{tikzcd}
      \Gamma(E_0) \ar[r, "\varphi"] & \Gamma(D) \ar[r, "\iota_*"] & \Gamma(V)
    \end{tikzcd}
  \end{equation*}
  is a morphism of modules, and thus by the Serre-Swan Theorem \ref{thm:gensingserreswan} is of the form $\Phi_*$ for a smooth morphism of VB-tepui $\Phi \colon E_0 \to V$. But such module morphisms are always smooth, and the map $\iota_*$ is an induction, therefore $\varphi$ is smooth, and $\varphi_*$ is well-defined.
  
  Moreover, for the same reason, $\varphi$ is continuous when one equips $\Gamma(E_0)$ with its standard Fréchet topology and $\Gamma(D)$ with the subspace topology of $\Gamma(V)$. This latter topology is also Fréchet, because $\Gamma(D)\subset \Gamma(V)$ is closed. In particular, $\varphi$ is a continuous bijection between Fréchet spaces, i.e., its inverse $\varphi^{-1}\colon \Gamma(D)\to \Gamma(E_0)$ is also continuous. This means that $\varphi^{-1}$ is diffeologically smooth, when equipping $\Gamma(E_0)$ and $\Gamma(D)$ with their $\sigma_\infty$-diffeologies. By Proposition \ref{prop:frechet-vs-diffeology}, the $\sigma_\infty$-diffeology of $\Gamma(E_0)$ coincides with its functional diffeology. Thus, it only remains to show that the $\sigma_\infty$-diffeology of $\Gamma(D)$ also coincides with its functional diffeology. for $f\colon U\to \Gamma(D)$ this can be seen as follows:
  \begin{align*}
  	&f\colon U\to \Gamma(D) &&\mathrm{is~ }\sigma_\infty\mathrm{-smooth}\\
  	\iff&i_*\circ f\colon U\to \Gamma(V) &&\mathrm{is~ }\sigma_\infty\mathrm{-smooth}\\
  	\iff&i_*\circ f\colon U\to \Gamma(V) &&\mathrm{is~ }\mathrm{smooth~w.r.t.~the ~functional~diffeology}\\
  	\iff&\widehat{i_*\circ f} \colon U\times M\to V &&\mathrm{is~ smooth}\\
  	\iff&\widehat{f}\colon U\times M\to D &&\mathrm{is~ smooth}\\
  	\iff&f\colon U\to \Gamma(D) &&\mathrm{is~ smooth~w.r.t.~the ~functional~diffeology.}  	
    \end{align*}
  Therefore, the pushforward by $\varphi$ thus gives an isomorphism of diffeological spaces
  \begin{equation*}
    \begin{tikzcd}
      \varphi_* \colon C^\infty(L, \Gamma(E_0)) \ar[r, "\cong"]&  C^\infty(L, \Gamma(D)),
    \end{tikzcd}
  \end{equation*}
and we conclude from Diagram \eqref{eq:diagramsurge} that $\alpha_D$ is surjective, as desired. Moreover, if $\sigma_1,\dots,\sigma_K$ are generators of $D$, then $\pi^*\sigma_1,\dots,\pi^*\sigma_k$ generate $\Gamma(\pi^*D)$.

The second important special case is when $\pi\colon P=U\to M$ is the inclusion of an open subset. In this case, we show surjectivity directly: let $f\in\Gamma(\pi^*D)=\Gamma(D|_U)$. There exists a smooth function $\chi\in C^\infty(M)$ which does not vanish on $U$, such that $\chi\cdot f$ can be extended (by zero) to a smooth function $\tilde f\in \Gamma(D)$ (cf. e.g. \cite[Chapter V, Lemma 6.1]{tougeronIdeauxFonctionsDifferentiables1972}). Then $(\frac{1}{\chi}\otimes \tilde f)\in C^\infty(U)\otimes_{C^\infty(M)}\Gamma(D)$ is a pre-image of $f$ along $\alpha_D$.

Now we address the general case, where $\pi\colon P \to M$ is a submersion. Let $\sigma_1,\dots ,\sigma_k$ be generators of $\Gamma(D)$. We claim that $\pi^*\sigma_1,\dots,\pi^*\sigma_k$ generate $\Gamma(\pi^*D)$. Locally, this is true, since any submersion is locally a projection over an open subset $U$ of $M$, and we know the $\pi^*(\sigma_i|_U)$ generate $\Gamma(\pi^*(D|_U))$ from the first and second step. The global statement follows by a partition of unity argument. But $\pi^*\sigma_k=\alpha_D(1\otimes \sigma_k)$, hence $\alpha_D$ is surjective onto generators, hence surjective.
\end{proof}

\begin{rk}
  At this point, for a vector bundle $V$ and a singular subbundle $D$, we have concluded that several possibly distinct smooth structures on the section spaces $\Gamma(V)$, $\Gamma(V/D)$, and $\Gamma(D)$, separately coincide. We gather these findings here. First, Proposition \ref{prop:frechet-vs-diffeology} ensures that the usual Fr\'{e}chet structure on $\Gamma(V)$ agrees with its canonical functional diffeology. This diffeology is also trivially the coefficient diffeology described in Remark \ref{rk:coefficient-diffeology}. Second, for $\Gamma(V/D)$, we have Lemma \ref{lem:lifting-sections} and Proposition \ref{prop:frechetmodfromsections}, where we deduce that the functional diffeology on $\Gamma(V/D)$ is the quotient diffeology it inherits from $\Gamma(V/D) \cong \Gamma(V)/\Gamma(D)$, and this diffeology agrees with the Fr\'{e}chet structure also inherited from $\Gamma(V/D) \cong \Gamma(V)/\Gamma(D)$. Once again, this is also the coefficient diffeology, as per Remark \ref{rk:coefficient-diffeology}. Finally, for $\Gamma(D)$, its functional diffeology is the subset diffeology it inherits from $\Gamma(D) \subseteq \Gamma(V)$, and in the proof of Proposition \ref{prop:welldefinepullback} above we saw that this agrees with the Fr\'{e}chet structure also inherited from $\Gamma(D) \subseteq \Gamma(V)$. In the case where $\Gamma(D)$ is locally finitely generated, this diffeology is also the coefficient diffeology, a fact in this case witnessed by the \emph{a fortiori} diffeomorphism $\Gamma(E_0) \cong \Gamma(D)$ from our Serre-Swan theorem \ref{thm:gensingserreswan}.  
\end{rk}
If $\pi$ is not a submersion, then it may be that $\gs(\pi^*E_0)$ is not equal to $\gs(\pi^*D)$, because the latter may not be finitely generated.

\begin{ex}
	Consider the trivial vector bundle $V \coloneqq \R \times \R \to \R$, with singular subbundle $D$ the image of the bundle endomorphism $(x,v) \mapsto (x,xv)$. Thus, the fiber of $D$ is $\R$ except at the origin, where it is trivial. The module $\gs(D)$ is projective. For $E_0$ the vector bundle such that $\gs(E_0)\cong\gs(D)$, one has that $\gs(\pi^*E_0)$ is finitely generated for any $\pi$. Nevertheless if we choose a function $\pi\colon \KR\fto\KR$ that vanishes on $(-\infty,0]$ and is positive otherwise (therefore not a submersion), then $\gs(\pi^*D)$ is not finitely generated, therefore $\gs(\pi^*E_0)$ is not isomorphic to  $\gs(\pi^*D)$.
\end{ex}

\begin{ex}
  One case where the finite generation of $\Gamma(D)$ is automatic is when $M$ is compact and $D$ is locally real-analytic (i.e., near any point of $M$ there are local coordinates and a trivialization of $V$ such that $D$ is generated by real-analytic sections). The proof of this is analogous to \cite[Proposition 2.3 and Theorem 2.4]{lavaulaurent-gengoux-el-al}.  For example, we may take the trivial vector bundle $V \coloneqq TS^2 = S^1 \times \R \to S^1$, with singular subbundle $D \coloneqq \mathrm{span}\{X_\theta \mid \theta \in S^1\}$, where $X$ is a vector field on $S^1$ that vanishes only at $1 \in S^1$. Then $D$ is locally real-analytic, and so $\Gamma(D)$ is finitely generated, and by Proposition \ref{prop:welldefinepullback} the pullbacks $\Gamma(\pi^*(V/D))$ and $\pi^*\Gamma(V/D)$ agree for any submersion $\pi\colon P \to S^1$.
\end{ex}

\section{Singular foliations and Lie tepui algebroids}
\label{sec:algarbroids}

In this last section, we turn to our original motivation for introducing tepui fibrations, namely applications in the theory of singular foliations. We begin with the definition of a tepui algebroid.

\subsection{Tepui algebroids}
\begin{defn} A \define{(Lie) tepui algebroid} is a triple $(A \to M, \rho, [\cdot,\cdot])$ where:
	\begin{itemize}
		\item $A \to M$ is VB-tepui;
		\item $\rho\colon A\to TM$ is a VB-tepui morphism over $1_M$. We call $\rho$ the \define{anchor};
		\item $[\cdot,\cdot]\colon \Gamma(A)\times \Gamma(A)\to \Gamma(A)$ is a Lie bracket (i.e., bilinear, skew-symmetric and satisfies the Jacobi identity);
	\end{itemize}	
	such that the following \define{Leibniz identity} holds:
	\begin{equation*}
		\text{for all }a,b\in\Gamma(A),\text{ and for all }f\in C^\infty(M), \text{ we have }[a,fb]=f[a,b]+\rho(a)(f)b
	\end{equation*}
      \end{defn}

\begin{ex}[Lie algebroids] Any Lie algebroid is a tepui algebroid. Indeed, the definition of tepui algebroid directly generalizes that of a Lie algebroid.
\end{ex}

\begin{rk}\label{rk:lierinhartsheaf} In view of the singular Serre-Swan Theorem \ref{thm:equiv1}, we can give a completely algebraic characterization of tepui algebroids: a tepui algebroid of globally bounded rank over $M$ is a fiber-determined finitely generated Lie-Rinehart algebra over $C^\infty(M)$. An analogous description is possible without the condition of globally bounded rank, using Theorem \ref{thm:gensingserreswan} and locally finitely generated sheaves of Lie-Rinehart algebras (or global locally finitely generated fiber-determined Lie Rinehart algebras). Such sheaves of Lie-Rinehart algebras were studied in \cite{villa}.
\end{rk}

Our main sources of examples for tepui algebroids are singular foliations. We recall the definition here.
\begin{defn}
	\label{def:singular-foliation}
	A \emph{singular foliation} on a smooth manifold, in the sense of \cite{AS09}, is a $C^\infty(M)$-submodule $\mathcal F$ of the module of compactly supported vector field $\mathfrak{X}_c(M)$ which is:
	\begin{itemize}
		\item involutive, i.e.,  $[\mathcal F,\mathcal F]\subset \mathcal F$, and
		\item locally finitely generated, i.e.,  each point of $M$ admits a neighbourhood $U$, such that $\mathcal F|_U$ is finitely generated as a module over $C^\infty(U)$. Here $\cl{F}|_U$ is the submodule of $\mathfrak{X}_c(U)$ consisting of those vector fields which are restrictions to $U$ of elements of $\cl{F}$ whose support is contained in $U$.
	\end{itemize}
\end{defn}
Singular foliations are naturally induced by Lie algebra actions on manifolds, Poisson structures, Lie algebroids, and much more. From a singular foliation, we may get a tepui algebroid.

\begin{ex}[Tepui algebroids from singular foliations]\label{ex:tepuialgfol}
Let $\mathcal F$ be a singular foliation. Then there is a VB-tepui $A = \operatorname{HolAlg}(\cl{F})$ such that $\Gamma(A)\cong \mathcal F$. To obtain it, we apply Theorem \ref{thm:gensingserreswan} to the global hull of $\mathcal F$.  The inclusion $\mathcal F\to\mathfrak X(M)$ then induces a morphism $\rho\colon A\to TM$. The corresponding map $\rho_*\colon\Gamma(A)\to \mathfrak X(M)$ is injective (with image the global hull of $\mathcal F$). In particular, the Lie bracket of $\mathfrak{X}(M)$ induces a bracket on $\Gamma(A)$ that turns $(A, \rho, [\cdot,\cdot])$ into a tepui algebroid. We call $A$ the \define{holonomy algebroid} of $\mathcal F$.
\end{ex}

\begin{rk} As a diffeological space, $\operatorname{HolAlg}(\cl{F})$ coincides with the ``fiberspace'' introduced in \cite[Definition 2.11]{GV22}. What we add here is the fact that it is a tepui fibration and a tepui algebroid globally.
\end{rk}

\begin{rk} Given a singular foliation $\mathcal F$ over $M$, one can construct its \emph{adiabatic foliation} $\mathcal T\mathcal F$ over $M\times \mathbb R$ \cite{androulidakisskanalisDNC}. In \cite[Proposition 3.6]{androulidakisskanalisDNC}, it is proven that the restriction of the holonomy groupoid of $\mathcal T\mathcal F$ to $M\times\{0\}$ is $\operatorname{HolAlg}(\cl{F})$. Together with Example \ref{ex:holonomy-groupoid}, this result also implies that $\operatorname{HolAlg}(\cl{F})\to M$ is a tepui fibration.
\end{rk}

Singular subalgebroids, introduced in \cite{singsub}, are a simultaneous generalization of wide subalgebroids and singular foliations. These also generate tepui algebroids.

\begin{ex}[Singular subalgebroids]
  \label{ex:singsub}
  Given a Lie algebroid $A \to M$, a \define{singular subalgebroid} is a locally finitely generated $C^\infty(M)$-submodule $\mathcal B$ of $\Gamma_c(A)$, which is stable under the Lie bracket. In complete analogy to the case of singular foliations, there is a tepui algebroid $\tilde{A}$ and a tepui algebroid morphism $A \to \tilde{A}$ such that $\Gamma_c(\tilde{A}) \cong \mathcal B$.
\end{ex}

For more examples, we appeal to Lemma \ref{lem:smoothsingvbtotepuis}; just as VB-tepui can arise as quotients of vector bundles, tepui algebroids can arise as quotients of (almost)-Lie algebroids.

\begin{defn}
  An \define{almost-Lie algebroid} is a triple $(A,\rho,[\cdot,\cdot])$, where:
  \begin{itemize}
  \item $A \to M$ is a vector bundle;
  \item $\rho\colon A\to TM$ is a VB-morphism over $1_M$;
  \item $[\cdot,\cdot]\colon \Gamma(A)\times \Gamma(A)\to \Gamma(A)$ is bilinear, skew-symmetric, and satisfies the Leibniz identity, and the following weakened Jacobi identity:
    	\begin{equation*}
		\textrm{for all } a,b,c \in \Gamma(A), \quad \operatorname{Jac}(a,b,c)\coloneqq [a,[b,c]]-[[a,b],c] - [b,[a,c]]\in \ker(\rho_*).
	\end{equation*}	
  \end{itemize}
\end{defn}

For details on almost-Lie algebroids, see \cite{LaurentGengouxLouisRyvkin}. Given an almost-Lie algebroid $A \to M$, we denote by $D_{\ker \rho_*}$ the smooth singular subbundle of $A$ spanned by $\ker(\rho_*)$. Note that the kernel of the bundle map $\rho$ contains $D_{\ker \rho_*}$, but they are not generally equal.

\begin{ex}
	\label{ex:almost-lie-quotient}
	Let $A \to M$ be an almost-Lie algebroid, and let $D$ be any smooth singular subbundle of $A$ containing $D_{\ker \rho_*}$, and such that $[\Gamma(D), \Gamma(A)] \subseteq \Gamma(D)$. We obtain a natural tepui algebroid structure on $A/ D$.
\end{ex}

In fact, any almost-Lie algebroid gives rise to a singular foliation, namely $\rho_*(\Gamma_c(A))$, and $A/ D_{\ker \rho_*}$ recovers the tepui algebroid $\operatorname{HolAlg}(\rho_*(\Gamma_c(A)))$ from Example \ref{ex:tepuialgfol}.

\begin{ex}[$L_\infty$-algebroids]\label{ex:linfty} An $L_\infty$-algebroid is given by a negatively graded vector bundle $E=\bigoplus_{i\in\mathbb N} E_{-i}$ over $M$, a VB-morphism $\rho\colon E_{-1}\to TM$, and graded-symmetric operators $l_i\colon \prod^i \Gamma(E)\to \Gamma(E)$ turning $\Gamma(E)$ into an $L_\infty$-algebra and satisfying:
	\begin{itemize}
		\item $l_i$ are $C^\infty(M)$-multilinear for $i\neq 2$;
		\item $l_2$ satisfies the Leibniz rule with respect to $\rho$ (where $\rho$ is extended by zero from $E_{-1}$ to $E$).
	\end{itemize}
	In particular, we obtain a complex of vector bundles
	\begin{equation*}
		\cdots \overset{l_1}\to E_{-k}\overset{l_1}\to \cdots \overset{l_1}\to E_{-2}\to E_{-1} \overset{\rho}\to TM
	\end{equation*}
	and a binary bracket $[\cdot,\cdot]\colon \Gamma(E_{-1})\times \Gamma(E_{-1})\to \Gamma(E_{-1})$ satisfying the Jacobi identity up to elements in $l_1(\Gamma(E_{-2}))$. This means that $E_{-1}/l_1(E_{-2})$ carries the structure of a tepui algebroid. If we apply this construction to the universal $L_\infty$-algebroid of a singular foliation (cf. \cite{lavaulaurent-gengoux-el-al}), we again obtain the holonomy algebroid $A(\mathcal F)$.
\end{ex}

In Example \ref{ex:almost-lie-quotient}, we saw that quotients of almost-Lie algebroids result in tepui algebroids. The converse also holds, in direct analogy to the case of VB-tepui in Proposition \ref{prop:tepui.is.quotient}.

\begin{lem}\label{lem:alwaysalmostexists}Let $(A,\rho, [\cdot,\cdot])$ be a tepui algebroid of globally bounded rank. Let $V$ be a vector bundle with singular smooth subbundle $D$ such that $A \cong V/D$, and let $\phi\colon V \to A$ be the quotient map. There exists an almost Lie algebroid structure $(V,\rho'=\rho\circ \phi,[\cdot,\cdot]')$ such that $\phi$ intertwines $[\cdot,\cdot]'$ with $[\cdot,\cdot]$. In particular, $A$ can be realized as the quotient of the almost Lie algebroid $V$ by the smooth singular subbundle $D=\ker(\phi)\subset V$.
\end{lem}
\begin{proof} This proof is very similar to the proof of the fact that any anchored bundle admits a compatible almost Lie bracket \cite[Proposition 2.2.4]{LaurentGengouxLouisRyvkin}.
	First, we find a complement $\tilde V$ to $V$, meaning a vector bundle $\tilde V \to M$ such that $V \oplus \tilde V = M\times \mathbb R^N$ for some $N\in\mathbb N$. Let $e_1,\dots,e_N$ be the standard frame of $M\times \mathbb R^N$. We set $\bar \phi \coloneqq \phi\circ \pi_V\colon M\times \mathbb R^N\to A$, where $\pi_V\colon M\times \mathbb R^N\to V$ is the projection. Since $\bar \phi_*\colon \Gamma(M\times \mathbb R^N)\to \Gamma(A)$ is surjective, we can find smooth functions $\tilde c_{ij}^k$ such that $[\bar\phi_*e_i,\bar\phi_*e_j]=\sum_k\tilde c_{ij}^k\bar\phi_*e_k$. We now set $c_{ij}^k=\frac{1}{2}\left(\tilde c_{ij}^k-\tilde c_{ji}^k\right)$ and define an almost-Lie bracket $[\cdot,\cdot]''$ on $M\times \mathbb R^N$ as the $\mathbb R$-bilinear map satisfying:
	\begin{equation*}
		[f_ie_i,f_je_j]'' \coloneqq f_if_j\sum_kc_{ij}^k e_k + f_i((\rho\circ \bar\phi)_*e_i)(f_j)e_j - f_j((\rho\circ \bar\phi)_*e_j)(f_i)e_i.
	\end{equation*}
	Finally, we set the bracket on $V$ to be:
	\begin{equation*}
		[a,b]'\coloneqq (\pi_V)_*[(i_V)_*a,(i_V)_*b]''
	\end{equation*}
	where $i_V\colon V\to M\times \mathbb R^N$ is the inclusion. By construction, this bracket is skew-symmetric and satisfies the Leibniz and weakened Jacobi rules. It is also compatible with the bracket on $A$.
\end{proof}

\begin{rk}
	In Lemma \ref{lem:alwaysalmostexists}, the sections of the smooth singular subbundle $\ker \phi$ form a Lie ideal in $\Gamma(V)$. Even when $V$ is a Lie algebroid, these are slightly more general than the ideals usually considered, since we do not require $\ker\phi$ to be a vector subbundle of $V$.
\end{rk}

In the definition of tepui algebroid, we have not explicitly required $\rho$ to be a Lie algebra homomorphism. As in the case of Lie algebroids, this is automatic.
\begin{lem}
	Let $A$ be a tepui algebroid. Then $\rho_*\colon \Gamma(A)\to \mathfrak{X}(M)$ is a Lie algebra homomorphism.
\end{lem}

\begin{proof} By Lemma \ref{lem:generic-regularity}, there is an open dense subset $U\subset M$, restricted to which $A$ has locally constant rank. Any point in $U$ admits a neighbourhood on which $A$ is a Lie algebroid, and in particular the restriction of $\rho$ is a Lie algebroid morphism, by \cite[Section 6.1]{Kosmann-Schwarzbach-Magri}. This means that for any $a,b\in\Gamma(A)$, the vector field $\rho_*[a,b]-[\rho_*a,\rho_*b]$ vanishes on an open and dense subset of $M$, and hence vanishes everywhere.
\end{proof}

As an immediate consequence, we obtain:
\begin{cor} Let $A$ be a tepui algebroid, then $\operatorname{Fol}(A) \coloneqq \rho(\Gamma_c(A))$ is a singular foliation.
\end{cor}

Together with Example \ref{ex:tepuialgfol}, this corollary implies that the category of tepui algebroids can be seen as an extension of the category of singular foliations. We explain this in the following remark.
\begin{rk}
	\label{thm:sing-fol-tep-alg-equiv}
	Singular foliations over a fixed base manifold $M$ form a category with inclusions as morphisms. Thus, a morphism $\mathcal F\to \mathcal F'$ exists if and only if $\mathcal F\subset \mathcal F'$, and it is unique in this case. On the other side, a morphism of tepui algebroids $A \to A'$ is a smooth VB-tepui map over the identity $1_M$, which intertwines the brackets and anchors of $A$ and $A'$.
	Let
        \begin{itemize}
        \item $\operatorname{HolAlg}$ denote the functor which maps a singular foliation to its holonomy algebroid;
        \item $\operatorname{Fol}$ denote the functor which maps a tepui algebroid to the singular foliation $\rho(\Gamma_c(A))$.
        \end{itemize}
Then $\operatorname{Fol}$ is left-adjoint to $\operatorname{HolAlg}$, i.e.,
	\begin{equation*}
		\operatorname{Hom}(\operatorname{Fol}(A),\mathcal F)=\operatorname{Hom}(A,\operatorname{HolAlg}(\cl{F})).
	\end{equation*}
	Moreover, $\operatorname{Fol} \circ \operatorname{HolAlg}$ is the identity.
\end{rk}

A singular foliation $\mathcal F$ induces a decomposition of $M$ into weakly embedded submanifolds called leaves, which by \cite{hermann1962} can be described in two equivalent ways:

\begin{itemize}
	\item The leaf through $p\in M$ of a singular foliation is the set of all points  $q\in M$, which can be reached by following a finite number of flows of elements in $\mathcal F$. 
	\item The leaf through $p\in M$ is the maximal connected weakly embedded submanifold $L\subset M$ such that $T_qL=\{X(q)~|~X\in\mathcal F\}\subset T_qM$ for any $q\in L$.
\end{itemize}

Let us recall that an (almost-)Lie algebroid $A\to M$ is called transitive when its anchor map $A\to TM$ is surjective. Similarly, we call a tepui algebroid \define{transitive} if its anchor map is surjective. If $A\to M$ is an almost-Lie algebroid, and $L$ is a leaf of its underlying singular foliation, then the restriction $A|_L$ defines a transitive almost-Lie algebroid.  Due to Lemma \ref{lem:alwaysalmostexists}, this implies:

\begin{lem} Let $A\to M$ be a tepui algebroid and $L$ a leaf of its underlying singular foliation. Then $A|_L$ is a transitive tepui algebroid.
\end{lem}

Transitive tepui algebroids can still be singular, as the following example shows.

\begin{ex}\label{ex:notlongsmooth}
	Let $A \coloneqq TM \times \R \to M$ be the Lie algebroid with anchor map $(v,t) \mapsto v$ and with Lie bracket on $\Gamma(A) \cong \mathfrak{X}(M) \oplus C^\infty(M,\R)$ given by
	\begin{equation}
		\label{eq:twisted-lie-bracket}
		[X\oplus f, Y \oplus g]_A \coloneqq [X,Y] \oplus (X(g) - Y(f)).
	\end{equation}
 	Fix some non-empty open subset $U \subseteq M$. Consider the smooth subbundle $D$ of $A$ given by
	\begin{equation*}
		D_m \coloneqq
		\begin{cases}
			\{0_m\} \times \{0\} &\text{if } m \in \overline{U} \\
			\{0_m\} \times \R &\text{otherwise}.
		\end{cases}
	\end{equation*}
	This is a smooth subbundle of $A$, because we may always find some smooth function $M \to \R$ which vanishes precisely on the closed set $\overline{U}$. Furthermore, the space of sections $\Gamma(D)$ is an ideal in $\Gamma(A)$. The condition $X \oplus f \in \Gamma(D)$ is equivalent to $X = 0$ and $f|_{\overline{U}} = 0$. Therefore, the Lie bracket of such $X \oplus f$ with $Y \oplus g$ becomes
	\begin{equation*}
		[0,Y] \oplus (0(g) - Y(f)) = 0 \oplus (-Y(f)).
	\end{equation*}
	Since $f$ vanishes on the open subset $U$, so does $Y(f)$. By continuity, $Y(f)$ also vanishes on $\overline{U}$, which is precisely the condition for $0 \oplus (-Y(f)) \in \Gamma(D)$. We conclude that $A/D$ is a transitive tepui algebroid that is not a Lie algebroid.
\end{ex}

\begin{rk}
  Example \ref{ex:notlongsmooth} provides a counterexample to the preprint \cite[Theorem 6.6]{villa} version [v2] of 2021, where it is stated that any locally finitely generated fiber-determined Lie-Rinehart sheaf over $M$ is \emph{adjoint integrable}. Adjoint integrability is defined in \cite[Definition 6.2]{villa}. Following Remark \ref{rk:lierinhartsheaf}, such a sheaf can be realized as the sheaf of sections of a tepui algebroid. In particular, the sheaf of sections of the tepui algebroid in Example \ref{ex:notlongsmooth} should be adjoint integrable, meaning that for any $a,b_0\in\Gamma(A)$ such that the flow of $\rho_*(a)$ is complete, the equation
\begin{align}\label{eq:adjint}
b'(t)=[b(t),a]
\end{align}
has a unique solution $b\colon \mathbb R\to \Gamma(A)$. Let us specify $M=\mathbb R$, with coordinate $x$, $U = (0,\infty)$, and $a=(\partial_x,0)$. Let $b_0$ be the zero section. Let $u$ be any function with support in $(-\infty,0)$. Then $b(t)(x)\coloneqq u(t+x)$ defines a solution of Equation \eqref{eq:adjint}, but this solution is not unique, showing that $\Gamma(A)$ is not adjoint integrable. We are grateful to Joel Villatoro for helpful discussions regarding this issue and confirming our observation. He has indicated that this issue will be addressed in a future version of the preprint.
\end{rk}

Since transitivity does not imply regularity in general, we single out those tepui algebroids which are non-singular along their leaves.

\begin{defn} We call a tepui algebroid $A$ \define{longitudinally smooth} if for any leaf $L$ of its underlying singular foliation, $A|_L$ is of constant rank, and hence is a Lie algebroid.
\end{defn}

Most tepui algebroids of interest to us are longitudinally smooth.

\begin{prop}\label{prop:longsmooth}
 Tepui algebroids coming from $L_\infty$-algebroids (as in Example \ref{ex:linfty}) or singular subalgebroids (as in Example \ref{ex:singsub}), and in particular those arising from singular foliations, are longitudinally smooth.
\end{prop}
\begin{proof}
For $L_\infty$-algebroids, this is a consequence of \cite[Lemma 1.9]{laurent-gengouxHolonomySingularLeaf2022}, and for singular subalgebroids, the proof is analogous to that of \cite[Proposition 2.11]{invitation-part2}.
\end{proof}

\begin{rk} There are certain similarities between the proofs for $L_\infty$-algebroids and singular subalgebroids. Both use the bracket on the almost-Lie algebroid $V$ that desingularizes the tepui algebroid $A$ to obtain a linear vector field, and then use its flow. For this procedure to work, the flow must preserve the kernel of the projection $\Gamma(V)\to \Gamma(A)$. This is true for different reasons:
	\begin{itemize}
		\item For $L_\infty$-algebroids this is a consequence of properties of the graded flow \cite[Lemma 1.9]{laurent-gengouxHolonomySingularLeaf2022}.
		\item For a singular subalgebroid $A\subset B$ of a Lie algebroid $B$, this follows from the fact that the kernel of $\Gamma(V)\to\Gamma(A)$ is also the kernel of a map $\Gamma(V)\to \Gamma(B)$ and the flow on $V$ can be related to a flow on $B$.
		\end{itemize}
	\end{rk}

\begin{ex} Let $\mathcal F$ be a singular foliation and $A=\operatorname{HolAlg}(\mathcal F)$ be the tepui algebroid from Example \ref{ex:tepuialgfol}, satisfying $\Gamma_c(A)=\mathcal F$. Then by Proposition \ref{prop:longsmooth}, for any leaf $L$ of $\cl{F}$ the bundle of vector spaces $A|_L$ is a Lie algebroid. The Lie algebroid $A|_L$ first appeared in \cite[Remark 1.16]{AS09} as the holonomy Lie algebroid of a leaf of the foliation. In \cite{debordLongitudinalSmoothnessHolonomy2013} it was shown that $A|_L$ is always integrable to a Lie groupoid, and that the holonomy groupoid of $\mathcal F$ restricted to $L$ provides this integration. In a forthcoming article \cite{tepuiholonomy}, we will develop the Lie theory for tepui fibrations and show how this $\operatorname{HolAlg}(\cl{F})$ is the differentiation of the holonomy groupoid from Example \ref{ex:holonomy-groupoid}.
\end{ex}

\subsection{On pullbacks}
In this subsection, we lay the groundwork to extend the original constructions of pullbacks for Lie algebroids in \cite{MK} to tepui algebroids. This presents some novel difficulties that are not present in the case of Lie algebroids. On the level of sections, our constructions recover the corresponding notions in the category of sheaves of Lie-Rinehart algebras handled in \cite{villa}.

Let $(A, \rho, [\cdot,\cdot])$ be a tepui algebroid, let $P$ be a manifold, and let $p\colon P\fto M$ be a smooth map. Consider the map of $C^\infty(P)$-modules
\begin{equation*}
 p^*\rho\colon \gi(P)\otimes_{\gi(M)}\gs(A) \fto \gs(p^* TM), \quad \Sigma_i g_i\otimes Y_i \mapsto \Sigma_i g_ip^*\rho(Y_i) 
\end{equation*}
and the $C^\infty(P)$-module $p^\bk\Gamma(A)$ given by the fiber product

\begin{equation*}
  \begin{tikzcd}
    p^\bk\gs(A)  \ar[r] \ar[d] & \gi(P)\otimes_{\gi(M)}\gs(A) \ar[d, "p^*\rho"] \\
    \gx(P) \ar[r, "dp"] & \gs(p^*TM ).
  \end{tikzcd}
\end{equation*}

As for Lie algebroids \cite{higginsAlgebraicConstructionsCategory1990}, there is a Lie bracket on $p^\bk\Gamma(A)$, given by:
\begin{equation*}
  \left[\, (X,\Sigma_{i} g_i\otimes Y_i)\,\, ,\,\, (X',\Sigma_{j} g_j'\otimes Y_j')\, \right]=\left(\, [X,X']\,\, , \, \,\Sigma_{i,j} g_i g_j' \otimes[Y_i,Y_j']+ \Sigma_{j} Xg_j'\otimes Y_j'- \Sigma_{i} X'g_i\otimes Y_i\, \right).
\end{equation*}

The Lie algebra $p^\bk\Gamma(A)$ can be used as a starting point for two important constructions that transfer the tepui algebroid $A$ to one over the space $P$: the pullback (or base change) and the action algebroid construction.

\begin{itemize}
\item It might happen that $p^\bk A \coloneqq TP \fiber{}{TM} A$ is a VB-tepui over $P$, and also that $\Gamma(p^\bk A) \cong p^\bk \Gamma(A)$. Then $(p^\bk A, \pr_1, [\cdot,\cdot])$ is a tepui algebroid over $P$, and is one candidate for a pullback. In the case that $A$ is a Lie algebroid, a sufficient --- but not necessary --- condition for $p^\bk A$ to be a well-defined Lie algebroid is that $p$ is a surjective submersion.
\item It might happen that we have an \define{infinitesimal action} of $A$ on $P$. Precisely, such an action is given by a map
  \begin{equation*}
    (-)^\dagger\colon \gs(A)\fto \gx(P),\quad  X\mapsto X^\dagger
  \end{equation*}
  satisfying the following four conditions:
  \begin{alignat*}{2}
    (X+Y)^{\dagger} &= X^{\dagger}+Y^{\dagger} & (uX)^\dagger &= (p^*u) X^\dagger \\
  [X,Y]^\dagger &= [X^\dagger, Y^\dagger] \quad & X^\dagger &\text{ is }p\text{-related to } \rho(X).   
  \end{alignat*}
  Then there is an injective map
  \begin{equation*}
    \gi(P)\otimes_{\gi(M)} \gs(A) \fto p^{!!}\gs(A), \quad \Sigma_i g_i \otimes Y_i \mapsto \left(\Sigma_i g_i Y_i^\dagger , \Sigma_i g_i \otimes Y_i  \right). 
  \end{equation*}
  Its image is a sub-Lie algebra of $p^\bk\Gamma(A)$, thus induces a bracket on $\gi(P)\otimes_{\gi(M)} \gs(A)$. Then, if $\Gamma(p^*A) \cong \gi(P)\otimes_{\gi(M)} \gs(A)$, we may equip $p^*A$ with the structure of a tepui algebroid. This is the case if $A$ is a Lie algebroid.
\end{itemize}

When trying to form either $p^\bk A$ or $p^*A$ for tepui algebroids, we are obstructed by the fact that $\gi(P)\otimes_{\gi(M)} \gs(A)$ and $\gs(p^*A)$ do not generally coincide.

\begin{ex}
\label{ex:horribleaaaaaarrrrhhhhhhhh} Let us consider the smooth manifold $M=\mathbb R$ with open subset $U=\{x>0\}$, and the tepui algebroid $A$ described in Example \ref{ex:notlongsmooth}. We now consider $P=\mathbb R$ with the smooth map $p\colon P\to M$ given by $p(y)\coloneqq y^2$, as in Example \ref{ex:algpullbnotgeom}. Then we can calculate that
$p^*A=p^*TM\oplus \tilde E $, where $\tilde E_y$ is trivial for $y\neq 0$, and $\mathbb R$ at $y=0$.  We thus have $p^{!!}A=TP\oplus \tilde E$. We claim that there can be no tepui bracket on $p^{!!}A$ which is compatible with $p$ (i.e., having the projection to $TP$ as the anchor). To see this, assume such a bracket exists. Let $\sigma$ be any non-zero section of $\tilde E$. Then $y\sigma=0$, but $[\partial_y,y\sigma]=\sigma+y[\partial_y,\sigma]\neq 0$, which is a contradiction.

For another perspective on the same phenomenon, observe that the invisible elements of $p^{!!}\Gamma(A)$ do not generally form a Lie ideal. In this example, if we pick $\tilde \sigma$ to be any element in $\Gamma(\ker \rho)\subset \Gamma(A)$ with non-trivial value at zero, we see that the invisible element $(0,y\otimes (0,\tilde \sigma))$ has non-trivial bracket with $(\partial_y, 2y\otimes (\partial_x,0))$.
\end{ex}

\begin{ex} Let $P=\mathbb R$, $M=\mathbb R^2$ and $p\colon P\to M$ be given by $p(t)\coloneqq (t,t)$. Consider the regular foliation $F\subset TM$ spanned by $\partial_x+\partial_y$ on $M$. We have an action of $F$ on $P$, given by restricting a vector field to the diagonal in $\R^2$. Let now $U \coloneqq [0,\infty)\times (-\infty, 0]\subset M$, and $E \to M$ be the VB-tepui which has trivial fibers over $M \smallsetminus U$ and $\mathbb R$-fibers over $U$, as constructed in Example \ref{ex:notlongsmooth}. We consider the tepui algebroid $A \coloneqq F\oplus E$ on $M$. This algebroid retains a natural transitive action on $P$. However, $p^*A$ can never be a tepui algebroid with bracket compatible with that on $A$, for the same reason as in the previous example.
\end{ex}

Having seen examples where pullback constructions fail, we turn to cases where we can perform the pullback. First, we have a result for some singular foliations.

\begin{prop} Let $\CF$ be a finitely generated singular foliation on a manifold $M$ and $p\colon P\fto M$ a submersion. If $A\fto M$ is the tepui algebroid of Example \ref{ex:tepuialgfol}, then $\gs(p^{!!}A)$ has a canonical Lie bracket, implying that the invisible elements form an ideal.
\end{prop}

\begin{proof}
	This uses the method found in \cite[Lemma 1.15]{GZ19} to prove that $\gs(p^{!!}A)$ is the VB-tepui of the pullback foliation $p^{-1}\CF$, and \cite[Proposition 1.10]{AS09} to prove that $p^{-1}\CF$ is a singular foliation.
	
	We will only check that $\rho(\gs(p^{!!}A))=p^{-1}\CF$, by double inclusion.
        \begin{itemize}
        \item For ``$\subset$,'' let $(X,Y)\subset\gs(p^{!!}A)=\gx(P)\times_{\gs(p^*TM)} \gs(p^*A)$. By Proposition \ref{prop:basechange} $Y$ is of the form
          \begin{equation*}
            Y=\sum f_i p^*Y_i \text{ where } f_i\in \gi(P) \y Y_i\in \gs(A).
          \end{equation*}
	Then since $p$ is a submersion, there exists sections $X_i$ such that $dp(X_i)=p^*\rho(Y_i)$, and so,
        \begin{equation*}
          X=\sum f_i X_i\in p^{-1}\CF
        \end{equation*}
        \item 	For the other inclusion ``$\supset$,'' let $X\in p^{-1}\CF$. Then 
          \begin{equation*}
            dp(X)=\sum f_i p^*\rho(Y_i) \text{ where } f_i\in \gi(P) \y Y_i\in \gs(A)
          \end{equation*}	
	and therefore $(X,\sum f_i p^*\rho(Y_i) )\in \gs(p^{!!}A)$.

        \end{itemize}

\end{proof}

In order to construct pullbacks in more general situations, we use the following lemma.

\begin{lem} The invisible elements of $p^{!!}\gs(A)$ are inside the ideal 
  \begin{equation*}
    I=\{0\}\times_{\gs(p^*TM)} \left(\gi(P)\otimes_{\gi(M)} \gs(A)\right).
  \end{equation*}
	 Therefore, if $\gi(P)\otimes_{\gi(M)} \gs(A)$ has no invisible elements, neither does $p^{!!}\gs(A)$.
\end{lem}
\begin{proof}
	The module $\gx(P)$ has no invisible elements, therefore the image of $p^{!!}\gs(A)\fto \gx(P)$ must be zero.
\end{proof}

This implies that we should focus on the invisible elements of $\gi(P)\otimes_{\gi(M)} \gs(A)$. For this, we use the following, derived from Proposition \ref{prop:welldefinepullback}.

\begin{cor}
\label{prop:pullback-exists} Let $V \to M$ be an almost-Lie algebroid with smooth singular subbundle $D$ containing $D_{\ker \rho_*}$, and denote by $A \coloneqq V/D$ the resulting tepui algebroid. Let $\pi\colon P \to M$ be a submersion, and suppose that $\Gamma(D)$ is projective. If there is an infinitesimal action of $\gs(A)$ on $P$, then $p^*A\fto P$ is a tepui algebroid.

Moreover if $\pi$ is transverse to the image of $\rho\colon A\fto TM$ then 
\begin{equation*}
  p^{!!}\gs(A)=\gs(p^{!!}A),
\end{equation*}
and $p^{!!}A\fto P$ is a tepui algebroid.
\end{cor}
\begin{proof}
	By Proposition \ref{prop:welldefinepullback} we have $\gi(P)\otimes_{\gi(M)}\gs(A)\cong \gs(p^*A)$ which immediately gives the first part.
	
	For the second part, note that:
        \begin{equation*}
          p^{!!}A = TP\times_{TM} A \cong TP\times_{TM} \frac{V}{D} \cong \frac{TP\times_{TM} V}{\{0\}\times D}
        \end{equation*}
	If $p$ is a submersion, then $TP\times_{TM} V$ is a vector bundle over $P$, and therefore $TP\times_{TM} A$ is a VB-tepui. Then, repeating the first part of the argument, we conclude that $p^{!!}\gs(A)=\gs(p^{!!}A)$.
      \end{proof}

\subsection{Singular foliations as tepui algebroids}
In this final subsection, we will discuss the advantages of interpreting singular foliations as tepui algebroids. Most of the results we obtain are not new, but we hope the reader finds their presentation in the language of tepui algebroids elegant.

\begin{lem}\label{lem:mortofol} Let $\mathcal F$ be a singular foliation and $A=\operatorname{HolAlg}(\mathcal F)\to M$ be its tepui algebroid as in Example \ref{ex:tepuialgfol}. Then:
\begin{itemize}
	\item The global hull $\hat{\mathcal F}$ of $\mathcal F$ is a Fréchet space. The foliation $\mathcal F$ itself is an LF-space.
	\item A time-dependent vector field in $\mathcal F$ in the sense of \cite[Definition 1.3.2]{LaurentGengouxLouisRyvkin} is simply a diffeologically smooth map $I\to \hat{\mathcal F}$. 	
\end{itemize}
\end{lem}
\begin{proof} 
The first item is a direct application of Proposition \ref{prop:frechetmodfromsections}. The LF-topology on $\Gamma(A)$ can be obtained as an inductive limit of the Fréchet topologies on $\Gamma_K(A)$, where $\Gamma_K(A)$ is the Fréchet space of sections of $A$ which have support in the compact subset $K\subset M$. The second item is immediate.
\end{proof}

We note that the Fréchet topology on $\hat{\mathcal F}$ (respectively the LF-topology on $\mathcal F$) need not be their topology as subspaces of $\mathfrak{X}(M)$ (respectively $\mathfrak{X}_c(M)$). In particular, a smooth path $X\colon I\to \mathfrak{X}(M)$ with $X_t\in\mathcal F$ for all $t$, might not constitute a time-dependent vector field of $\mathcal F$ in the sense of \cite[Definition 1.3.2]{LaurentGengouxLouisRyvkin}. For details, see \cite[Remark 1.3.3]{LaurentGengouxLouisRyvkin}. When  $\hat{\mathcal F}$ is a closed subspace of $\mathfrak X(M)$, then its subspace topology coincides with the topology induced by $\Gamma(A)$, as a consequence of the closed graph theorem.

\begin{lem} \label{lem:prodfol} Let $\mathcal F_1,\mathcal F_2$ be singular foliations and $A_i \coloneqq \operatorname{HolAlg}(\cl{F}_i)$ be their holonomy algebroids. The product singular foliation of $\mathcal F_1$ and $\mathcal F_2$ on $M_1\times M_2$  is uniquely described by the product tepui algebroid $A_1\times A_2\to M_1\times M_2$.
\end{lem}

The normal bundle of a singular foliation can also be interpreted as a VB-tepui.

\begin{ex}\label{ex:transversemod} Let $\mathcal F$ be a singular foliation. We call the VB-tepui $\mathcal T \coloneqq TM/D_{\mathcal F}$, where $D_{\mathcal F}$ is the smooth singular subbundle spanned by $\mathcal F$, the \define{normal bundle} to $\mathcal F$.
\end{ex}

\begin{defn} Given any singular foliation $\CF$  with normal bundle ${\mathcal T}$, we define the \define{Bott connection} $\nabla$ using the formula 
\[\nabla_X (W+\CF) =[X, W]+\CF,\]
for any $X\in \CF$ and $W+\CF \in \Gamma({\mathcal T})$.
\end{defn}
The Bott connection satisfies the usual conditions for a connection; this is, for any $f\in \gi(M)$, $X\in \CF$, and $W+\CF\in \Gamma(\mathcal T)$,
\begin{align*}
  \nabla_{fX} (W+\CF) &= f [X, W]+\CF, \\
  \nabla_{X} (fW+\CF) &= df(X) W + f [X, W]+\CF.
\end{align*}

Because both $\CF$ and $\Gamma({\mathcal T})$ have the coefficient diffeology, the map
\begin{equation*}
 \nabla\colon \CF\times \gs({\mathcal T})\fto \gs({\mathcal T}),\quad  (X,W+\CF)\mapsto \nabla_X (W+\CF) 
\end{equation*}
is smooth. To be more precise, let $\bar{Z}\colon \KR^n\fto \CF\times \gs({\mathcal T})$ be a smooth map, so that $\bar{Z}=(\bar{X},\bar{W})$ where $\bar{X}\colon \KR^n\fto \CF$ and $\bar{W}\colon \KR^n\fto \gs({\mathcal T})$ are smooth maps. Write:
\begin{align*}
  \bar{X} &= \sum_i f_i X_i  \text{ for } f_i\in \gi(\KR^n\times M) \text{ and } X_i \in \CF \\
  \bar{W} &= \sum_j g_j W_j +\CF  \text{ for } g_i\in \gi(\KR^n\times M) \text{ and } W_i+\CF \in \gs({\mathcal T}).
\end{align*}
Then 
  \begin{equation*}
    \nabla_{\bar X} \bar{W}= \sum_{i,j} f_i \left( dg_j(X_i) W_j + g_j [X_i,W_j]\right)+\CF,
  \end{equation*}
which is again smooth in the coefficient diffeology of $\gs({\mathcal T})$.
        
Moreover, parallel transport is possible in this setting. Using \cite[Lemma 2.3.1]{GV22} and \cite{gayu},
we know that for any curve $[0,1]\fto \CF\times M$ with $t\mapsto (X_t,\gamma(t))$, satisfying $X_t(\gamma(t))=\gamma'(t)$ (these are called $\CF$-paths in \cite{GV22}) we have $[X_t,\CF]\subset \CF$, and the rank of $\CF_{\gamma(t)}$ is constant in $t$. Then $\gamma^*A(\CF)$ is a vector bundle over $[0,1]$. This implies that the rank of ${\mathcal T}_{\gamma(t)}$ is constant in $t$, and $\gamma^*{\mathcal T}$ is a vector bundle over $[0,1]$. Then for any $W_0\in {\mathcal T}_{\gamma(0)}$ it is possible to solve the differential equation
          \begin{equation*}
            \nabla_{X_t} W(t)=0; \quad  W(0)=W_0
          \end{equation*}
	using the same techniques as for parallel transport on regular vector bundles. This gives a unique element $W_1=W(1)\in {\mathcal T}_{\gamma(1)}$, which we interpret at the transport of $W$ along $(X,\gamma)$.

\bibliographystyle{alpha}

\bibliography{alga.bib}

\end{document}